\documentclass{amsart}

\usepackage{amsmath}
\usepackage{amsfonts}

\newcommand{\Z}{\mathbb{Z}}
\newcommand{\F}{\mathbb{F}}

\newcommand{\Dh}{\mathcal{D}_H}

\newcommand{\U}[2]{\underbrace{#1, \dots, #1}_{#2}}
\newcommand{\n}[0]{\text{$\mathbf{n}$}}
\newcommand{\Sl}[0]{\mathrm{sl}}

\DeclareMathOperator{\ad}{ad}
\DeclareMathOperator{\Der}{Der}
\DeclareMathOperator{\SkDer}{SkDer}

\newtheorem{dummy}{Dummy}

\numberwithin{dummy}{section}
\numberwithin{equation}{section}

\newtheorem{lemma}[dummy]{Lemma}
\newtheorem{theorem}[dummy]{Theorem}

\theoremstyle{definition}
\newtheorem{definition}[dummy]{Definition}

\theoremstyle{remark}
\newtheorem{rem}[dummy]{Remark}



\begin{document}

\bibliographystyle{amsalpha}

\title[Gradings of non-graded Hamiltonian Lie algebras]%
{Gradings of non-graded\\
  Hamiltonian Lie algebras\\}

\author{A.~Caranti}

\address{%
  Dipartimento di Matematica\\
  Universit\`a degli Studi di Trento\\
  via Sommarive 14\\
  I-38050 Povo (Trento)\\
  Italy}

\email{caranti@science.unitn.it}

\author{S.~Mattarei}

\address{%
  Dipartimento di Matematica\\
  Universit\`a degli Studi di Trento\\
  via Sommarive 14\\
  I-38050 Povo (Trento)\\
  Italy}

\email{{mattarei@science.unitn.it}}

\thanks{The  authors  are grateful  to  Ministero dell'Istruzione, dell'Universit\`a  e
  della  Ricerca, Italy,  for  financial  support to  the
  project ``Graded Lie algebras  and pro-$p$-groups of finite width''.}

\maketitle

\begin{abstract}
A thin Lie algebra is a Lie algebra graded over the positive integers satisfying a certain narrowness condition.
We describe several cyclic grading of the modular Hamiltonian Lie algebras $H(2\colon\n;\omega_2)$
(of dimension one less than a power of $p$)
from which we construct infinite-dimensional thin Lie algebras.
In the process we provide an explicit identification of $H(2\colon\n;\omega_2)$ with a Block
algebra. We also compute its second cohomology group and its derivation algebra (in arbitrary prime characteristic).
\end{abstract}



\thispagestyle{empty}

\section{Introduction}\label{sec:introduction}

In the last three decades the interest of researchers in finite $p$-groups has increasingly extended to pro-$p$ groups.
This trend was initiated by Leedham-Green and Newman in 1980, who proposed in~\cite{L-GN}
one way of getting around the universally believed impossibility of a classification of $p$-groups up to isomorphism.
One of their intuitions was that of using the {\em coclass} rather than the (nilpotency) class of  $p$-group
as a fundamental invariant.
Since the coclass of a group of order $p^n$ is defined as the difference between $n$ and the class of the group,
this change has no real effect unless one focuses on families of $p$-groups rather than single $p$-groups.
A special role is then played by pro-$p$ groups, to which the definition of coclass extends naturally,
and which represent entire families of finite $p$-groups as the sets of their finite quotients.
In particular, having finite coclass for a pro-$p$ group $G$ means that all lower central quotients
$\gamma_i(G)/\gamma_{i+1}(G)$ have order at most $p$ from some point on.
The five {\em coclass conjectures} advanced in~\cite{L-GN} are now theorems thanks to the efforts
of several authors, culminating in~\cite{L-G} and~\cite{Sha:coclass}.
They give information on pro-$p$ groups of finite coclass
(Conjecture~C, the simplest to state, claiming that every pro-$p$ group of finite coclass is soluble),
but also asymptotic information on families of finite $p$-groups of fixed coclass.

Lie-theoretic methods occupy an important place in the theory of $p$-groups,
and also in the proofs of the coclass conjectures.
The oldest and simplest such method is associating the graded Lie ring
$\bigoplus_{i=1}^\infty \gamma_i(G)/\gamma_{i+1}(G)$
with the lower central series of a (pro-) $p$-group.
In many interesting cases the Lie ring is actually a Lie algebra over the field of $p$ elements $\F_p$.
This approach was already present in disguise in the pioneering work of Blackburn
on $p$-groups of coclass one (better known as $p$-groups {\em of maximal class})
and the subsequent work of Leedham-Green and McKay (see~\cite{L-GMcK:maximal-classification} and the references therein),
which inspired the formulation of the coclass conjectures.
Lie rings or algebras
associated with pro-$p$ groups in this way are graded over the positive integers,
and are generated by their component of degree one.
The coclass conjectures have natural analogues for graded Lie algebras over an arbitrary field defined by these properties,
independently of the connection with pro-$p$ groups,
but in this new context all the conjectures turn out to be false.
This already occurs in the simplest case of coclass one, as
Shalev constructed  in~\cite{Sha:max} the first examples of  insoluble graded
Lie    algebras   of    maximal   class.
Here, by a graded Lie algebra of maximal class we mean a Lie
algebra $L$ which is graded over the positive integers,
\begin{equation}\label{eq:graded}
  L = \bigoplus_{k = 1}^{\infty} L_{k},
\end{equation}
where $L_{1}$ has dimension $2$, $L_{k}$ has dimension $1$ for $k >1$,
and $L_{k+1} = [L_{k}, L_{1}]$ for all $k \ge 1$.

Shalev's construction starts from certain
finite-dimensional simple modular  Lie algebras, originally introduced
by Albert and Frank~\cite{AF},
and applies a {\em loop  algebra}
construction (strictly speaking, taking the positive part of a twisted loop algebra,
see Subsection~\ref{subsec:maximal} for details).
As a consequence, the resulting graded Lie algebras of maximal class, despite being insoluble,
have a kind of periodic structure, in a precise sense.
Therefore, one could still hope that each of them is uniquely determined by a suitable
finite-dimensional quotient, as is the case with pro-$p$ groups of finite coclass
(because each $p$-group of finite coclass $r$ and class large enough, depending on $p$ and $r$,
is a quotient of a unique infinite pro-$p$ group of coclass $r$).
The investigation carried out in~\cite{CMN} showed that this is not the case.
In fact, it turned out that most graded Lie algebras of maximal class are not periodic.
Nevertheless, Shalev's algebras occupy a unique place in the description of the
graded Lie algebras of maximal class over fields of odd characteristic,
which have been completely classified in~\cite{CN}.
An analogous classification over fields of characteristic two will soon appear in~\cite{Ju:maximal}.

After the coclass conjectures for pro-$p$ groups
were proved, other invariants have been suggested, which may possibly lead to
a finer classification of $p$-groups than the rough one provided by the coclass theorems.
We do not need the precise definitions of these invariants, which are called {\em width,} {\em obliquity} and {\em rank},
and we refer the interested reader to Chapter~12 of the book~\cite{L-GMcKay} by Leedham-Green and McKay.
The simplest nontrivial case in terms of these invariants consists of the
pro-$p$ groups of width two and obliquity zero.
These groups, which do not have finite coclass in general,
have also been given the name {\em thin} in~\cite{Br,BCS} (originally in the finite case)
because of a narrowness characterization in terms of their lattice of (closed) normal subgroups.
Some of our terminology, like the {\em diamonds} introduced in the next paragraph,
originates from that point of view.

An approach to thin groups via the associated graded Lie algebra was taken in \cite{CMNS}.
The lower central factors in a thin group are elementary abelian of rank at most two.
Those of rank two, in the group or in its associated graded Lie algebra, are called {\em diamonds}.
There is of course a diamond $G/\gamma_2(G)$ on top of a thin group $G$, and if this is the only diamond
then $G$ has maximal class.
Otherwise, it follows from the theory of $p$-groups of maximal class that in a thin group the second diamond
occurs in class at most $p$.
In \cite{CMNS} we proved that the associated Lie algebra has bounded dimension,
except when the second diamond occurs in class (or degree) 3, 5 or $p$.
Each of these cases occurs for certain infinite pro-$p$ groups, both $p$-adic analytic and not.
For example, Sylow pro-$p$ subgroups of $\mathrm{SL}(2,\Z_p)$ or $\mathrm{SL}(2,\F_p[[t]])$,
or certain `nonsplit' versions of them, are thin pro-$p$ groups with second diamond in class 3.
A detailed description of these groups in the $p$-adic analytic cases is given in Section~12.2 of~\cite{L-GMcKay}.
Thin pro-$p$ groups with second diamond in class 5 can be realized similarly starting with
certain linear groups of type $A_2$ over local fields, see~\cite{KL-GP} or~\cite{Mat:thin-groups}.
Finally, the second diamond in class $p$ occurs for one of the wildest known pro-$p$ groups,
the Nottingham group (which is thin for $p>2$), described in Section~12.4 of~\cite{L-GMcKay}, for example.

A crucial fact for the investigation carried out in~\cite{CMNS}
was that the condition of a pro-$p$ group $G$ having obliquity zero
can be verified on the associated graded Lie ring, which in this case is a Lie algebra over $\F_p$.
In this context a more convenient formulation of the condition is the {\em covering property}.
A graded Lie algebra $L$ as in~\eqref{eq:graded}, over an arbitrary field and thus not necessarily associated with a group,
is called {\em thin} if $L_{1}$ has dimension $2$ and the
following \emph{covering property} holds:
\begin{equation}\label{eq:covering}
  \text{for all $k \ge 1$, and all $u \in L_{k}$, with $u \ne 0$, we
  have $L_{k+1} = [u, L_{1}]$.}
\end{equation}
It follows that $L$ is generated by $L_1$ as a Lie algebra
and that all homogeneous components have dimension at most two.
(See Definition~\ref{def:thin} and the comments that follow for more details.)
The arguments of \cite{CMNS} have been extended in~\cite{CaMa:thin,AviJur}
to show that the second diamond in an infinite-dimensional
thin Lie algebra (or one of finite dimension large enough)
can only occur in degree $3$, $5$, $q$ or $2q-1$, for some power $q$ of the characteristic $p$ of the underlying field.
It follows from~\cite{CMNS}
that there are, up to isomorphism and with the possible exception of very small characteristics,
one or two (depending on the ground field) infinite-dimensional thin Lie algebras
with second diamond in degree 3 and no diamond in degree 4,
and one with second diamond in degree 5.
Thus, they are the graded Lie algebras associated with some of the thin pro-$p$ groups mentioned earlier.
Machine computations have shown that each of the cases
where the second diamond is in degree $q$ or $2q-1$
splits into a number of subcases, depending on the location of the further diamonds and their {\em types}
(see Subsection~\ref{subsec:thin}).
Several of these subcases have been investigated in various papers, having in mind
a classification of all infinite-dimensional thin Lie algebras as a distant goal.
We refer to the paper~\cite{CaMa:Nottingham}, and to the references mentioned there,
for a discussion of thin Lie algebras with second diamond in degree $q$,
and to~\cite{CaMa:thin,AviMat:-1} for those with second diamond in degree $2q-1$.
Here we restrict ourselves to some general and informal comments on the type of results
which have been proved so far.

The results in this subject typically come in pairs, of rather different flavour:
namely, a uniqueness theorem and an existence theorem.
The former states that a certain initial structure of an infinite-dimensional thin Lie algebra
(formulated in terms of the location of the first few diamonds and their {\em types}),
determines the algebra completely (within the class of thin Lie algebras).
More precisely, a certain specified finite-dimensional thin Lie algebra is a quotient of a unique infinite-dimensional
thin Lie algebra $L$.
This is proved by producing a finite presentation for a central extension $M$
(broadly speaking, the universal central extension of $L$).
Usually $L$ itself is not finitely presented,
see Remark~\ref{rem:AFS-cocycles} for a specific instance of this phenomenon.
The existence result consists in the explicit construction of $L$,
as a loop algebra of some finite-dimensional Lie algebra $S$, with respect to a suitable cyclic grading,
and sometimes with the intervention of an outer derivation of $S$.

The latter type of result brings in an interesting connection with (usually simple) finite-dimensional modular Lie algebras,
certain cyclic gradings of them, and their derivations.
Their second cohomology group (with values in the trivial module) also plays a role, being closely related
with the centre of $M$.
Apart from the classical algebras of types $A_1$ and $A_2$ used in the construction of
thin Lie algebras with second diamond in degree $3$ and $5$ (see~\cite{CMNS,Car:thin_addendum}),
all the remaining cases involve (non-classical) simple modular Lie algebras of Cartan type,
namely Zassenhaus algebras and Hamiltonian algebras of various types.
We recall the definitions of these algebras in Sections~\ref{sec:Cartan} and~\ref{sec:Hamiltonian},
and point out in Remark~\ref{rem:field} other notations in use for them.
In particular, infinite-dimensional thin Lie algebras with second diamond in degree $q$
have been constructed as loop algebras of Zassenhaus algebras
(which have dimension a power of $p$, see~\cite{Car:Nottingham,Car:thin_addendum,Car:Zassenhaus-three}),
and Hamiltonian algebras of the types $H(2:\n;\omega_0)=H(2:\n)$
(the graded simple Hamiltonian algebras, of dimension two less than a power of $p$, see~\cite{Avi})
and $H(2:\n;\omega_1)$ (which are Albert-Zassenhaus algebras and have dimension
a power of $p$, see~\cite{AviMat:Nottingham}).

A preliminary version of the present paper, which predated and inspired some of the other papers cited here,
had as its main goal the construction of some
infinite-dimensional thin Lie algebras with second diamond in degree $2q-1$ as loop algebras of Hamiltonian algebras
$H(2:\n;\omega_2)$, which have dimension one less than a power of $p$.
(In fact, a construction for those thin Lie algebras had already been given in~\cite{CaMa:thin},
but as loop algebras of certain finite-dimensional Lie algebras defined `ad hoc'.)
The paper has somehow expanded after we have realised that some of our result may be of interest
independently of their application to thin Lie algebras.

Now we describe the contents of this paper in some detail.
Our results are presented in Sections~\ref{sec:cohomology}--\ref{sec:Andrea}.
Sections~\ref{sec:Cartan} and~\ref{sec:Hamiltonian} are expository
and include information on low characteristics which is not easily accessible in the literature.
Since much motivation for studying Lie algebras of
maximal class and related classes of `narrow' Lie algebras like thin Lie
algebras comes from analogous classes of
(pro-)$p$-groups, we have written the expository sections with the aim of making the paper accessible to the group
theorist with little knowledge of modular Lie algebras.

We have mentioned earlier the relevance of the second cohomology group (with values in the trivial module)
of a finite-dimensional Lie algebra with presentations of the corresponding loop algebra.
Therefore, we  determine in Section~\ref{sec:cohomology}
the second  cohomology  group of  the
algebras $H(2\colon\n;\omega_2)$.
We do this according to the classical method used in
Farnsteiner~\cite{Far:central}, which relates the group to the derivations of the algebra,
exploiting the presence of a nonsingular associative form.
The result is  surely well known to experts, but
we have been unable to find a suitable reference in the literature.
We have also collected in Section~\ref{sec:cohomology}
known information about derivations and the second cohomology group of
$H(2\colon\n;\omega_0)$, in all positive characteristics.

In Section~\ref{sec:gradings} we show how various cyclic gradings of $H(2:\n;\omega_2)$ 
can be constructed in a natural way, and how they relate to graded Lie algebras of maximal class
and to certain thin algebras closely connected to them, which were also studied in~\cite{CaMa:thin}.
Note  that $H(2:\n;\omega_2)$ has a natural filtration, the {\em standard} filtration,
inherited by the natural $\Z$-grading of the divided power algebra
on which it acts.
The filtration is not induced by a grading, however, in contrast to $H(2:\n;\omega_0)$,
which is a {\em graded} Lie algebra of Cartan type.
(Hence the cheap pun in the title.)

The first manifestation of a connection between `narrow' infinite-dimensional graded Lie algebras
and certain finite-dimensional simple Lie algebras was in Shalev's construction in~\cite{Sha:max} of insoluble
graded Lie algebras of maximal class, which we have mentioned earlier.
The finite-dimensional simple algebras of Albert and Frank used by Shalev belong to the larger class of
Block algebras, introduced in~\cite{Bl}.
It is known that Block algebras are Lie algebras of Hamiltonian Cartan type and, in particular,
that the algebras of Albert and Frank are Hamiltonian algebras $H(2\colon\n;\omega_2)$,
at least for $p>3$.
(We discuss this further and provide appropriate references towards the end of Section~\ref{sec:Hamiltonian}.)
In Theorem~\ref{thm:AF=Hamiltonian} we give an explicit
isomorphism of  the  Hamiltonian algebras  $H(2\colon\n;\omega_2)$ with
the algebras of Albert and Frank
which is valid in every positive characteristic.

The crucial property of the algebras of Albert and Frank exploited by Shalev in~\cite{Sha:max} is that
they admit a cyclic grading with one-dimensional components, and a nonsingular derivation
which permutes them transitively.
Benkart, Kostrikin and Kuznetsov proved in~\cite{BKK} and~\cite{KoK3} that the only finite-dimensional simple Lie algebras
with this property, over an algebraically closed field of characteristic $p>7$,
are the Hamiltonian algebras $H(2\colon\n;\omega_2)$.
The proof rests on the classification of simple modular Lie algebras,
which causes the restriction on $p$,
but now a classification-free proof of this result for $p>2$ is a consequence of~\cite{CN}.
A more thorough discussion, which includes the case of characteristic two,
will be found in Section~\ref{sec:AF}.
The cyclic grading of $H(2\colon\n;\omega_2)$ and the nonsingular derivation are also the ingredients
for discovering the isomorphism with a Block algebra exhibited in Theorem~\ref{thm:AF=Hamiltonian}.

The determination of the second cohomology group of $H(2\colon\n;\omega_2)$ in Section~\ref{sec:cohomology}
depends on a knowledge of its derivation algebra.
Since this piece of information is not easily available in low characteristics,
we fill this gap by exploiting the other incarnation of $H(2\colon\n;\omega_2)$,
as originally defined by Albert and Frank.
Thus, in Theorem~\ref{thm:AF-derivations} we compute the derivation algebra of the algebras of Albert and Frank
in arbitrary positive characteristic, by
suitably modifying Block's original proof~\cite{Bl},
which was valid for $p>3$ only (although for a larger class of algebras).

In~\cite{CaMa:thin} we constructed various thin Lie algebras with second diamond in degree $2q-1$.
One can assign a {\em type} to each diamond of such algebras, taking values in the underlying field plus $\infty$,
in such a way that the locations and types of the diamonds determine the isomorphism type of the algebra
(see Subsection~\ref{subsec:thin} for more details).
The thin Lie algebras with all diamonds of type $\infty$ turn out to be closely connected with graded
Lie algebras of maximal class.
In particular, there is such a thin Lie algebra associated with each of the graded Lie algebras
of maximal class constructed by Shalev in~\cite{Sha:max}.
We describe that in~Subsection~\ref{subsec:thin}
as a loop algebra of $H(2\colon\n;\omega_2)$ with respect to a suitable grading.
In~\cite{CaMa:thin} we also constructed thin Lie algebras with second diamond in degree $2q-1$
and all diamonds of finite types, as loop algebras of certain Block algebras.
The Block algebras used in~\cite{CaMa:thin} are actually isomorphic with algebras of Albert and Frank
(being simple Block algebras with $G=G_0$, see Section~\ref{sec:Hamiltonian},
and according to known results, for example~\cite[Lemma~1.8.3]{BlWil:rank-two}),
but were presented there in a different basis.
As we have mentioned earlier, the original motivation for the present paper was finding an explicit identification
of those algebras with Hamiltonian algebras $H(2:\n;\omega_2)$
and describing a corresponding cyclic grading of the latter.
We realize these two tasks in Sections~\ref{sec:Andrea} and~\ref{sec:thin}, respectively.
We should mention that coexistence of diamonds of both finite and infinite types is possible in the same thin Lie algebra.
Such algebras are constructed in~\cite{AviMat:-1}, as loop algebras of Hamiltonian algebras $H(2:\n;\omega_2)$
with respect to yet another grading.

We are grateful to the referees for their comments.

\section{Generalities about some Lie algebras of Cartan type}\label{sec:Cartan}

The Hamiltonian Lie algebras form one of the four families
of Lie algebras of Cartan type $W$, $S$, $H$, $K$.
Definitions of these families can be found
in the recent book of Strade~\cite{Strade:book}.
Here we will limit ourselves to a discussion of the
general Lie algebra of Cartan type $W(m:\n)$,
with a special attention for the Zassenhaus algebra $W(1:n)$, and of two types of Hamiltonian algebras
in the next section.

Let $\F$ be a field of prime characteristic $p$,
let $\n=(n_1,\ldots,n_m)$ be an $m$-tuple
of positive integers, and put $n=n_1+\cdots+n_m$.
The algebra of $m$ divided powers truncated at $\n$,
denoted by $\F[x_1,\ldots,x_m;n_1,\ldots,n_m]$
or $\F[m:\n]$ for brevity,
is the $\F$-vector space of formal $\F$-linear combinations of monomials
$x_1^{(i_1)}\cdots x_m^{(i_m)}$ with $0\le i_j<p^{n_j}$,
with multiplication defined by
$x_j^{(k)} x_j^{(l)}=\binom{k+l}{k} x_j^{(k+l)}$,
and extended by linearity and by postulating commutativity and
associativity of the multiplication.

Note that, as an algebra, $\F[m:\n]$ is determined up to
isomorphism by its dimension $p^n$, where $n=n_1+\cdots+n_m$. In
fact, it coincides (up to notation) with the free associative and
commutative algebra on the generators $x_j^{(p^{k_j})}$, for
$0\le j\le m$, $0\le k_j<n_j$, subject to the law $x^p=0$ (that is, with
$\F[x_1,\ldots,x_n]/(x_1^p,\ldots,x_n^p)$); this is easily seen by
using Lucas' theorem \cite{Lucas} to compute the binomial
coefficient $\binom{k+l}{k}$.
In particular, a derivation of $\F[m:\n]$ can be defined by
sending the given free generators to arbitrarily chosen elements of
$\F[m:\n]$ and extending by the Leibniz rule; furthermore, every derivation
is obtained in this way.

However, $\F[m:\n]$ comes equipped with an additional structure,
namely a set of {\em divided power maps}, which tie the various
$p$-(divided) powers of the same variable together~\cite{Kac}.
We will not need to know any detail about the divided power maps, except
that the definition of {\em special} derivations given in \cite{Kac}
or \cite{SF} singles out exactly
those derivations of $\F[m:\n]$ which are compatible with
the divided power maps in a natural sense.
It turns out that the special derivations of $\F[m:\n]$ are
those of the form
$D=f_1\, \partial/\partial x_1+\cdots+f_m\, \partial/\partial x_m$
with $f_j\in\F[m:\n]$, thus acting as
$Dx_j^{(k)}=f_jx_j^{(k-1)}$ (where $x_j^{(-1)}=0$).

The Lie algebra of special derivations of $\F[m:\n]$
(which coincides with the full derivation algebra of $\F[m:\n]$
only when $n_1=\cdots=n_m=1$)
is denoted by $W(m:\n)$ and is called the
{\em general Lie algebra of Cartan type}
(or {\em generalized Jacobson-Witt algebra}).
It is simple of dimension $mp^n$, unless $m=1$ and $p=2$.
A grading of $W(m:\n)$ over $\Z^m$ is inherited from the natural
$\Z^m$-grading of $\F[m:\n]$, but it is the grading of $\F[m:\n]$
given by total degree of monomials which induces the most important grading
of $W(m:\n)$, called the {\em standard grading:}
$W(m:\n)=\bigoplus_{i=-1}^{r} L_i$, where
$L_i$ is the subspace spanned by the derivations
$x_1^{(i_1)}\cdots x_m^{(i_m)}\, \partial/\partial x_j$
with $i_1+\cdots+i_m=i+1$ and $j=1,\ldots,m$
(where $r=p^{n_1}+\cdots+p^{n_m}-m-1$).

\begin{rem}
What we have described is the unique generalized
Jacobson-Witt algebra, for fixed $m$ and $\n$, provided the field $\F$ is
algebraically closed; since we have taken $\F$
arbitrary of prime characteristic, $W(m:\n)$ is just one of possibly
many $\F$-forms of the generalized Jacobson-Witt algebra,
see~\cite{Wat}, for example.
A similar proviso applies to the Hamiltonian algebras which we will
describe in the next section, see~\cite{SerWil}
for a determination of the forms in the restricted case.
\end{rem}

In this paper we will actually only need the Zassenhaus algebras
$W(1:n)$, as they occur as distinguished subalgebras
of the Hamiltonian algebras which we will consider.
In this case the components of the standard grading are one-dimensional,
$L_i$ being generated by $E_i= x^{(i+1)}\, d/dx$, for $i=-1,\ldots,r$
(where $r=p^n-2$).
Direct computation shows that
\begin{equation*}
[E_i,E_j]=
\left(
\binom{i+j+1}{j}-
\binom{i+j+1}{i}
\right)
E_{i+j}.
\end{equation*}
In particular, we have $[E_{-1},E_j]=E_{j-1}$,
and $[E_0,E_j]=jE_j$.

The Zassenhaus algebra has also a grading over
(the additive group of) $\F_{p^n}$, with graded basis
consisting of the elements $e_{\alpha}$, for $\alpha\in\F_{p^n}$,
which satisfy
\begin{equation*}
[e_{\alpha},e_{\beta}]=
(\beta-\alpha)
e_{\alpha+\beta}.
\end{equation*}
In particular, note that $[e_0,e_{\alpha}]=\alpha e_{\alpha}$.
The bases $\{E_i\}$ and $\{e_{\alpha}\}$ of $W(1;n)$ are sometimes referred to in the literature
(at least when $n=1$) as a {\em proper basis} and a {\em group basis}, respectively.
One way to obtain the group basis from the proper basis is noting that
$E_{-1}+E_r$ spans a Cartan subalgebra of $W(1;n)$, and computing
the corresponding Cartan decomposition.
Since $\ad (E_{-1}+E_r)$ permutes
$E_{r-1},E_{r-2},\ldots E_1,E_0,E_{-1}+2E_r$
cyclically, one quickly finds the formulas
\begin{equation*}
\begin{cases}
  e_0= E_{-1}+E_r
  &\\
  e_{\alpha}= E_r+ \sum_{i=-1}^{r} \alpha^{i+1}E_i & \text{for
    $\alpha\in\F_{p^n}^{\ast}$},
\end{cases}
\end{equation*}
where $r=p^n-2$.
Note that the first formula is a special case of the
second formula if we stipulate that $0^0=1$.

The inverse formulas are
\begin{equation*}
\begin{cases}
  E_{-1}= e_0+ \sum_{\alpha} e_{\alpha}
  &\\
  E_{i}= -\sum_{\alpha} \alpha^{r-i} e_{\alpha} & \text{for
    $i=0,\ldots,r$},
\end{cases}
\end{equation*}
where the summations are for $\alpha\in\F_{p^n}$, and again for the
case $i=r$ we understand $\alpha^0=1$ for any $\alpha\in\F_{p^n}$,
whence $E_{r}= -\sum_{\alpha} e_{\alpha}$.

We note in passing that the Zassenhaus algebra is not simple when $p=2$,
but has a unique non-trivial ideal, namely
$\langle E_i\mid i\not=r\rangle=
\langle e_{\alpha}\mid\alpha\not=0\rangle$,
which we will refer to as {\em the simple Zassenhaus algebra}.
The above transition formulas are clearly valid in this case, too.
As we will point out in the next section, a simple Zassenhaus algebra in characteristic two can also be regarded
as a Hamiltonian algebra (and, in turn, as a Block algebra).

Since the problem of inverting formulas similar to those which relate
the bases $\{e_{\alpha}\}$ and $\{E_i\}$ of the Zassenhaus algebra
will occur repeatedly in this paper, we record the solution explicitly.
To simplify notation and computations it will be useful
to set $0^0=1$ once and for all.
The customary rules $\alpha^i\beta^i=(\alpha\beta)^i$ and
$\alpha^i\alpha^j=\alpha^{i+j}$ now hold for $\alpha,\beta$ in a field and
$i,j$ non-negative integers.
Note that with this convention the expression
$\sum_{\alpha\in \F_q} \alpha^{j}$,
where $\F_q$ is the finite field of $q$ elements,
becomes meaningful for every non-negative integer $j$;
its value is $-1$ if $j$ is a positive multiple of $q-1$,
and $0$ otherwise.

\begin{lemma}
The linear relations
\begin{equation*}
a_{\alpha}=\sum_{j=0}^{q-1} \alpha^{j} b_j,
\quad
\text{for $\alpha\in \F_q$,}
\end{equation*}
between elements
$a_{\alpha}$ $(\alpha\in \F_q)$ and
$b_j$ $(j=0,\ldots,q-1)$
of any vector space over $\F_q$, 
are equivalent to the relations
\begin{equation*}
\begin{cases}
  b_0=a_0, \\
  b_j=-\sum_{\alpha\in \F_q} \alpha^{q-1-j} a_{\alpha}, \quad
  \text{for $j=1,\ldots,q-1$.}
\end{cases}
\end{equation*}
\end{lemma}

\begin{proof}
If $\omega$ is a primitive $n$-th root of unity in any
field, then the sets of formulas
\begin{equation*}
a_i=\sum_{j=1}^n \omega^{ij} b_j,
\quad\text{and}\quad
b_j=\frac{1}{n}\sum_{i=1}^n \omega^{-ij} a_i,
\end{equation*}
relating subsets $\{a_i\mid i=1,\ldots,n\}$ and
$\{b_j\mid j=1,\ldots,n\}$ of any vector space over that field,
are inverse of each other.
This is an instance of a Fourier transform and its inverse over a cyclic
group of order $n$, and can be easily proved using the fact that
$\sum_{j=1}^n \omega^{ij}$ equals $n$ if $i$ is a multiple of $n$,
and $0$ otherwise.

In particular, taking as $\omega$ a generator of $\F_q^{\ast}$,
we obtain that the sets of formulas
\begin{equation*}
a_{\alpha}=\sum_{j=1}^{q-1} \alpha^{j} b_j,
\quad\text{and}\quad
b_j=-\sum_{\alpha\in \F_q^{\ast}}
\alpha^{-j} a_{\alpha},
\end{equation*}
which relate elements
$a_{\alpha}$ $(\alpha\in \F_q)$ and
$b_j$ $(j=0,\ldots,q-1)$
of any vector space over $\F_q$, are
inverse of each other.

Now consider the formulas
$a_{\alpha}=\sum_{j=0}^{q-1} \alpha^{j} b_j$
for $\alpha\in \F_q$,
note that one of them is $a_0=b_0$,
and write the remaining ones in the form
$a_{\alpha}-a_0=\sum_{j=1}^{q-1} \alpha^{j}b_j$.
As we have seen above, these formulas can be inverted, and yield
\begin{equation*}
b_j=
-\sum_{\alpha\in \F_q^{\ast}}
\alpha^{-j} (a_{\alpha}-a_0)=
-\delta_{j,q-1}\, a_0
-\sum_{\alpha\in \F_q^{\ast}}
\alpha^{-j} a_{\alpha}=
-\sum_{\alpha\in \F_q}
\alpha^{q-1-j} a_{\alpha},
\end{equation*}
for $j=1,\ldots,q-1$.
\end{proof}

\section{Hamiltonian algebras and Block algebras}\label{sec:Hamiltonian}

The Lie algebras of Cartan type $S$, $H$ and $K$ are defined as
subalgebras of the generalized Jacobson-Witt algebra $W(m:\n)$,
and depend on a choice of a certain differential form $\omega$
(or, equivalently, of a certain automorphism of $W(m:\n)$).
For this to make sense in general one must complete the algebra of
divided powers $\F[m:\n]$ to an algebra of divided power series.
However, this will not be necessary to define the only Hamiltonian
algebras which we will consider in this paper, namely
$H(2:\n;\omega_j)$ for $j=0$ or $2$.
Note that the Hamiltonian algebras in two variables can also be
considered as belonging to the Cartan series of special algebras
(and thus be denoted by $S(2:\n;\omega_j)$).
A rather condensed description of all four classes of Lie algebras of Cartan type,
but complete with all the relevant references,
can be found in~\cite{BKK}, to which we also conform our notation.
(See Remark~\ref{rem:field} concerning our notation.)
For a more extensive discussion see~\cite{Strade:book}.

As in the previous section, we assume only that $\F$ is a field of prime characteristic $p$,
and point out which statements need restrictions on $\F$ as we go along.
Let $\F[2:\n]=\F[x,y;n_1,n_2]$ be
the algebra of divided powers in two variables $x,y$
of heights $\n=(n_1,n_2)$.
It will be convenient to put $\bar x=x^{(p^{n_1}-1)}$,
$\bar y=y^{(p^{n_2}-1)}$, and $e=\bar x\bar y$.
Then $H(2:\n;\omega_j)$ can be defined as the second derived
algebra of
\[
\tilde H(2:\n;\omega_j)=
\{ D\in W(2:\n)\mid
D\omega_j=0 \},
\]
where
$\omega_0=dx\wedge dy$ and
$\omega_2=(1-e)\, dx\wedge dy$.
Note that the (formal) differential forms here are simply elements of the
exterior algebra on the set $\{dx,dy\}$ over $\F[2:\n]$.
In particular, the space of differential $2$-forms
is the free $\F[2:\n]$-module
on the basis $\{dx\wedge dy\}$, and is a $W(2:\n)$-module via
\[
D(f\,dx\wedge dy)=(Df)\,dx\wedge dy+f\,d(Dx)\wedge dy+f\,dx\wedge d(Dy)
\]
for $D\in W(2:\n)$,
where $df=(\partial f/\partial x)\,dx+(\partial f/\partial y)\,dy$.
When dealing with derivations of $H(2:\n;\omega_j)$ it will be useful
to consider the larger algebra
\[
C\tilde H(2:\n;\omega_j)=
\{ D\in W(2:\n)\mid
D\omega_j=c\omega_j,\, c\in\F \},
\]
which contains $\tilde H(2:\n;\omega_j)$ as an ideal of
codimension one and zero for $j=0,2$, respectively.

\begin{rem}\label{rem:field}
We should mention that the notation for the Lie algebras of Cartan type
is not uniform in the literature.
In particular, $H(-)$ sometimes denotes what we have indicated with
$\tilde H(-)$ here, and so one has to take the (first or) second derived algebra $H(-)^{(2)}$
to obtain the simple algebra.
Furthermore, special, Hamiltonian and contact algebras can also be obtained as subalgebras
of the generalized Jacobson-Witt algebra (on an algebra of divided power series, in general)
by means of certain automorphisms $\Phi$, instead
of differential forms $\omega$.
For example, our $H(2:\n;\omega_0)$ and $H(2:\n;\omega_2)$
are denoted by
$H(2;\underline{n})^{(2)}=H(2;\underline{n};\mathrm{id})^{(2)}$
and $H(2;\underline{n};\Phi(\tau))^{(1)}$ in the book~\cite{Strade:book},
and similarly in many papers.
We also note that, strictly speaking, the notation
$H(2:\n;\omega_i)$ which we use here
would be only justified when working over an algebraically closed field
(and of characteristic large enough).
This is because only in that case it can be shown that any form $\omega$ defining a Hamiltonian algebra
can be assumed to have certain specific forms $\omega_0$, $\omega_1$, $\omega_2$
(see~\cite{BKK} and the reference therein for the most general results,
but~\cite[Corollary~2]{Wil:type-S} suffices for the Hamiltonian algebras $H(2:\n;\omega_i)$
under consideration here). 
Whenever we consider $H(2:\n;\omega_i)$ over an arbitrary field in this paper,
we refer to the specific form defined above.
\end{rem}

Since the space of differential $2$-forms on $\F[2:\n]$ has a natural
structure of graded module for the $\Z^2$-graded Lie algebra
$W(2:\n)$, and $\omega_0$ is a homogeneous element with respect to this
grading, $\tilde H(2:\n;\omega_0)$ and $C\tilde H(2:\n;\omega_0)$
are graded subalgebras of $W(2:\n)$ with respect to the $\Z^2$-grading.
In particular, they are also graded subalgebras with respect
to the standard grading of $W(2:\n)$, and thus they acquire what is called
their {\em standard grading}.
They are usually referred to as the {\it graded}
Lie algebras of Hamiltonian type, in contrast to
their relatives with respect to forms of type $\omega_1$ and $\omega_2$,
which are only {\em filtered}.

Thus, in determining an explicit expression for the generic
element $D$ of these subalgebras one may assume that $D$ is homogeneous
with respect to the $\Z^2$-grading of $W(2:\n)$.
A simple computation as in~\cite[pp.~255--257]{KoSa},
or \cite[p.~162ff.]{SF} (where the assumption that $p>2$
is not used before Theorem~4.5) shows that
$\tilde H(2:\n;\omega_0)$ consists of
the derivations of $\F[2:\n]$ of the form
$\Dh(f)= f_y\, \partial/\partial x-f_x\, \partial/\partial y$
for some
$f\in\tilde P(2:\n;\omega_0)=
\F[2:\n]
\oplus\langle
x^{(p^{n_1})},\,y^{(p^{n_2})}
\rangle$
(where $f_x$ and $f_y$ stand for
$\partial f/\partial x$ and $\partial f/\partial y$, respectively),
and that
$C\tilde H(2:\n;\omega_0)=\tilde H(2:\n;\omega_0)
\oplus\langle x\,\partial /\partial x\rangle$.
The latter can be written in the more symmetric form $\tilde H(2:\n;\omega_0)
\oplus\langle x\,\partial /\partial x+y\,\partial /\partial y\rangle$
if $p>2$.
Since
\begin{equation*}
[\Dh(f),\Dh(g)]=
(f_yg_x-f_xg_y)_y\frac{\partial}{\partial x}-
(f_yg_x-f_xg_y)_x\frac{\partial}{\partial y}=
\Dh(\Dh(f)(g)),
\end{equation*}
the map
$\Dh$
is a homomorphism from the $\tilde H(2:\n;\omega_0)$-module
$\tilde P(2:\n;\omega_0)$
onto the adjoint module for $\tilde H(2:\n;\omega_0)$, with kernel $\langle 1\rangle$.
The associative (and commutative) algebra $\tilde P(2:\n;\omega_0)$
can be endowed with an additional structure of Lie algebra
with respect to the {\em Poisson bracket}
$\{f,g\}=\Dh(f)(g)=
f_yg_x-f_xg_y$,
and the map $\Dh$ yields a Lie algebra isomorphism from
$\tilde P(2:\n;\omega_0)/\langle 1\rangle$ onto $\tilde H(2:\n;\omega_0)$.
(Note that our notation for the map $\Dh$ and, consequently, for the Poisson bracket,
differs in sign from that of~\cite{SF}, and agrees with~\cite{Kos:beginnings}
or~\cite{BO} instead.)
Under this isomorphism, the second derived subalgebra
$H(2:\n;\omega_0)$ of $\tilde H(2:\n;\omega_0)$ corresponds to the subalgebra
\[
P(2:\n;\omega_0)/\langle 1\rangle=\langle
x^{(i)}y^{(j)}\in
\F[2:\n]\mid
x^{(i)}y^{(j)}\not= e
\rangle
/\langle 1\rangle,
\]
of dimension $p^n-2$, where $n=n_1+n_2$.
This is a simple Lie algebra if $p>2$,
see~\cite{KoSa} or~\cite{SF}. 
In characteristic two it is simple provided $n_1>1$ and $n_2>1$,
as one can prove along the lines of Theorem~3.5 or~Theorem~4.5 of~\cite[Chapter 4]{SF}.
However,
$P(2:(1,n_2);\omega_0)/\langle 1\rangle$ has
$\langle y^{(j)}\mid 0\le j<2^{n_2}\rangle/\langle 1\rangle$
as an ideal.
In fact, it is the split extension of a simple Zassenhaus algebra by its adjoint module.
It can also be viewed as the tensor product of a simple Zassenhaus algebra with
the algebra of divided powers $\F[z:1]$.

Since the Lie algebra homomorphism $\Dh$
is also a homomorphism of $\tilde H(2:\n;\omega_0)$-modules,
the action of $\ad \Dh(f)$ as an inner derivation of the Lie algebra structure of $\tilde P(2:\n;\omega_0)$
coincides with the action of $\Dh(f)$ as a derivation of the associative algebra structure of $\tilde P(2:\n;\omega_0)$.
This will be useful when computing with derivations of $H(2:\n;\omega_0)$
in Sections~\ref{sec:cohomology} and~\ref{sec:gradings}.
For this reason we will simply regard $\Dh(f)$ as a derivation
of the Lie algebra $\tilde H(2:\n;\omega_0)$
(rather than the more cumbersome notation $\ad \Dh(f)$).
A word of caution, however:
this does not extend to derivations of $\tilde P(2:\n;\omega_0)$
which are not inner.
In fact, $\Dh$ is not a homomorphism of $C\tilde H(2:\n;\omega_0)$-modules, because
$[D,\Dh(x^{(i)}y^{(j)})]=
(i+j-2)\Dh(x^{(i)}y^{(j)})$
while
$D(x^{(i)}y^{(j)})=
(i+j)x^{(i)}y^{(j)}$,
for
$D=x\,\partial/\partial x+y\,\partial/\partial y$.

In this paper we will find convenient to always talk about $\tilde H(2:\n;\omega_0)$
while actually carrying out explicit computations inside
$\tilde P(2:\n;\omega_0)/\langle 1\rangle$ with the Poisson bracket
(and similarly for $\tilde H(2:\n;\omega_2)$, later).
Writing $x^{(i)}y^{(j)}$ for $x^{(i)}y^{(j)}+\langle 1\rangle$ we have
\begin{eqnarray*}
\{x^{(i)}
y^{(j)},
x^{(k)}
y^{(l)}\}
&=&
x^{(i)}
y^{(j-1)}
x^{(k-1)}
y^{(l)}-
x^{(i-1)}
y^{(j)}
x^{(k)}
y^{(l-1)}
\\
&=&
N(i,j,k,l)\,
x^{(i+k-1)}
y^{(j+l-1)},
\end{eqnarray*}
where
\[
N(i,j,k,l)=
{i+k-1\choose i}
{j+l-1 \choose j-1}-
{i+k-1 \choose i-1}
{j+l-1 \choose j}.
\]

The $\Z^2$-grading of $L=\tilde H(2:\n;\omega_0)$ is
$L=\bigoplus_{(i,j)\in\Z^2}L_{i,j}$, where
$L_{i,j}=\langle x^{(i+1)}y^{(j+1)}\rangle$ in the Poisson bracket notation,
and the standard grading is
$L=\bigoplus_{k\in\Z}\bar L_k$, where $\bar L_k=\sum_{i+j=k}L_{i,j}$
consists of all homogeneous polynomials of degree
$k+2$.
Note that $\bar L_k$ is trivial unless
$-1\leq k\leq p^{n_1}+p^{n_2}-4$
(or unless
$-1\leq k\leq p^{n_1}+p^{n_2}-5$,
if we restrict our attention to
$H(2:\n;\omega_0)$).

We consider now $C\tilde H(2:\n;\omega_2)$.
Although the form $\omega_2$ is not homogeneous
with respect to the $\Z^2$-grading, it becomes so
with respect to the grading obtained by viewing it modulo $(p^{n_1}-1,p^{n_2}-1)$.
Thus, $C\tilde H(2:\n;\omega_2)$ is a graded subalgebra of $W(2:\n)$
with respect to its $A$-grading, where
$A=\Z^2/\langle(p^{n_1}-1,p^{n_2}-1)\rangle$.
We will examine some specializations of this grading
in Section~\ref{sec:gradings}.
Consideration of $A$-homogeneous elements makes it straightforward to determine
an explicit form for the elements of $C\tilde H(2:\n;\omega_2)$.
One finds that
$\tilde H(2:\n;\omega_2)$ coincides with $C\tilde H(2:\n;\omega_2)$,
and can be identified with
$\tilde P(2:\n;\omega_2)/\langle 1\rangle=
(\F[2:\n]
\oplus\langle
x^{(p^{n_1})},\,y^{(p^{n_2})}
\rangle)
/\langle 1\rangle
$
with the Poisson bracket
$
\{f,g\}=(1+e)
(f_yg_x-f_xg_y).
$
The derived subalgebra $H(2:\n;\omega_2)$ of
$\tilde H(2:\n;\omega_2)$ has dimension $p^n-1$ and corresponds to
$
\F[2:\n]
/\langle 1\rangle
$
with the Poisson bracket.
It is simple (in every characteristic), as will follow from its identification,
given in Section~\ref{sec:AF}, with a certain Block algebra, whose simplicity was proved in~\cite{Bl}.
In characteristic two $H(2:(1,n);\omega_2)$
is isomorphic with the simple Zassenhaus algebra of dimension $2^{n+1}-1$,
an isomorphism being obtained by mapping $xy^{(j)}\mapsto E_{j-1}$ and $y^{(j)}\mapsto E_{j+2^{n}-2}$.
A curious consequence of this isomorphism is that in characteristic two $H(2:(1,n);\omega_2)$
can be embedded in $H(2:(1,n+1);\omega_2)$ as a subalgebra, namely as the simple Zassenhaus subalgebra
$\langle xy^{(j)}:j=0,\ldots,p^{n+1}-2\rangle$.

The Poisson bracket of monomials
\begin{equation*}
\{x^{(i)}
y^{(j)},
x^{(k)}
y^{(l)}\}
=
(1+e)\, N(i,j,k,l)\,
x^{(i+k-1)}
y^{(j+l-1)}
\end{equation*}
for $\tilde H(2:\n;\omega_2)$
coincides with that for $\tilde H(2:\n;\omega_0)$
except for the products $\{y,x\}=-\{x,y\}=e$.
This shows that the $A$-grading of $\tilde H(2:\n;\omega_2)$
cannot be lifted to a $\Z^2$-grading.
For the same reason, $\tilde H(2:\n;\omega_2)$ is not $\Z$-graded
by the subspaces $\bar L_k=\sum_{i+j=k}L_{i,j}$ defined as before.
However, it is filtered by the subspaces
$\bar L^k=\sum_{h\geq k} \bar L_h$, that is,
$L=\bar L^{-1}\supseteq \bar L^0
\supseteq \bar L^1\supseteq\cdots$
(where $\bar L^k=0$ for
$k>{p^{n_1}+p^{n_2}-4}$), and
$[\bar L^h,\bar L^k]\subseteq \bar L^{h+k}$.
This is called the {\em standard filtration} of $\tilde H(2:\n;\omega_2)$,
and the graded Lie algebra associated with it is
$\tilde H(2:\n;\omega_0)$.

We should point out here that both Hamiltonian algebras
$H(2:\n;\omega_0)$ and $H(2:\n;\omega_2)$ were originally constructed in a different way.
In fact, after being introduced first in~\cite{AF} among other examples, they became special cases
of a more general construction due to Block~\cite{Bl}.
We briefly recall only a special case of Block's construction
which is relevant to the present paper,
and we refer to~\cite{Bl} or~\cite[p.~110]{Sel} for full generality.

Let $G$ be an elementary abelian $p$-group of order $p^n$, written additively,
let $\delta\in G$, and let $f:G\times G\to G$ be a non-singular biadditive function
of the form
$f(\alpha,\beta)=g(\alpha)\,h(\beta)-g(\beta)\,h(\alpha)$
for some additive functions $g,h:G\to G$.
A vector space $L$ over a field $F$ of characteristic $p$,
with basis $\{u_{\alpha}\mid \alpha\in G\}$ in bijective correspondence
with the elements of $G$, becomes a Lie algebra by defining a multiplication
on the basis elements via
$[u_{\alpha},u_{\beta}]=f(\alpha,\beta)\,u_{\alpha+\beta-\delta}$
and extending linearly.
The element $u_0$ is central in $L$, and the elements $u_\alpha$ with $\alpha\not=\delta$
span an ideal of $L$.
If $\delta=0$ the ideal
$\langle u_\alpha\mid\alpha\not=0\rangle$ is a simple Lie algebra,
and if $\delta\not=0$ the quotient
$\langle u_\alpha\mid\alpha\not=\delta\rangle/\langle u_0\rangle$ is simple.
In both cases the simple algebra is called a {\em Block algebra}.
(These special cases of Block's construction had already been introduced by Albert and Frank in~\cite{AF},
and denoted by $\mathcal{L}_0$ and $\mathcal{L}_{\delta}$ there;
in this paper we refer to the algebras $\mathcal{L}_{\delta}$ with $\delta\neq 0$ as {\em algebras of Albert and Frank},
conforming to~\cite{Sha:max,CMN,CN}.)

It is known that if $\F$ is algebraically closed of characteristic $p>3$
the above special cases of Block's construction yield exactly the Hamiltonian algebras
$H(2:\n;\omega_2)$ if $\delta=0$,
and the algebras
$H(2:\n;\omega_0)$ if $\delta\not=0$.
For example, this is stated in~\cite[Lemma~1.8.3]{BlWil:rank-two} under the blanket assumption
of that paper that $p>7$, but the proof given there is seen to be valid for $p>3$.
(In particular, one ingredient of that proof,
namely~\cite[Corollary~2]{Wil:type-S}, was originally proved for $p>5$;
however, it is now a special case of more general results in~\cite{BGOSW} or~\cite{Skr:Hamiltonian}
which assume only $p>3$.)

Note that the method of proof of~\cite[Lemma~1.8.3]{BlWil:rank-two}
does not easily produce explicit realizations of $H(2:\n;\omega_0)$ and $H(2:\n;\omega_2)$
(with respect to the given forms) as Block algebras.
In fact, in essence (using automorphisms $\Phi$ rather than forms $\omega$), it
shows that for appropriate choices of the form $\omega$
the Hamiltonian algebra $H(2:\n;\omega)$ is a Block algebra of dimension
$p^n-2$ or $p^n-1$, and then appeals to~\cite[Theorem~1.8.1]{BlWil:rank-two}
(which quotes~\cite[Corollary~2]{Wil:type-S})
to conclude that $H(2:\n;\omega)\cong H(2:\n;\omega_i)$ for $i=0$ or $2$, respectively.
In Section~\ref{sec:AF} of the present paper we do give an explicit realization of $H(2:\n;\omega_2)$
as a Block algebra, and we do that for arbitrary prime characteristic $p$ (thus including $2$ and $3$).
We mention that, more generally, it was announced in~\cite{KoK1}
and proved in~\cite{KoK3} that every algebra $H(m:\n;\omega_2)$ is a Block algebra.

\section{The second cohomology group of $H(2:\n;\omega_2)$}\label{sec:cohomology}

In this section we appeal to some results which were formulated under the assumption that the ground field
is algebraically closed, or at least perfect;
since derivations and cohomology are essentially independent of the ground field,
these assumptions are immaterial here in view of~Remark~\ref{rem:field}.
We assume that the ground field has odd characteristic.
At some stage in the discussion we also need to assume that the characteristic $p$ is greater than three
(see Remark~\ref{rem:char3} for the case $p=3$),
but our main result, Theorem~\ref{thm:cocycles}, does not depend on this assumption.
Finally, we deal with the case of characteristic two in Remarks~\ref{rem:char2-derivations}
and~\ref{rem:char2-cocycles}.

The dimensions of the second cohomology groups $H^2(L,\F)$ of some
graded Lie algebras of Cartan type with values in the trivial module
were computed in \cite{Far:central}
(but see also~\cite{Dzhu:central}).
In particular, according to
\cite[Theorem~2.4]{Far:central}, $H^2(L,\F)$ has dimension $n_1+n_2+1$ for
the graded Hamiltonian algebra $L=H(2:\n;\omega_0)$.
Here we compute $H^2(L,\F)$ for $L=H(2:\n;\omega_2)$ and show that
it has dimension $n_1+n_2$.
(In the special case $\n=(1,1)$ this can essentially be found in~\cite[Theorem~6.3]{Str:Block}.)

Following~\cite{Far:central}, we briefly recall the classical method (see \cite[p.~102]{Sel})
employed there to compute $H^2(L,\F)$
from the space of outer derivations of $L$, in presence of a nondegenerate
associative form on $L$.
In addition, we exhibit a basis of
$H^2(L,\F)$ for $L=H(2:\n;\omega_0)$.
According to \cite[Proposition~1.3]{Far:central}, for any Lie algebra $L$
over a field $\F$ there is an injective homomorphism
$H^2(L,\F)\to H^1(L,L^{\ast})$,
where $L^{\ast}$ denotes the dual of the adjoint module of $L$.
This monomorphism is induced by the map
$\varphi\mapsto D_{\varphi}$ which sends a $2$-cocycle $\varphi\in
Z^2(L,\F)$ to the derivation $D_{\varphi}:L\to L^{\ast}$ with
$D_{\varphi}(\xi)=\varphi(\xi,\cdot)$.
Furthermore, the image of the monomorphism consists
of the cohomology classes represented by {\it skew} derivations,
that is to say, derivations $D:L\to L^{\ast}$ which satisfy
$D(\xi)(\eta)=-D(\eta)(\xi)$, for all $\xi,\eta\in L$.

Now assume that $L$ possesses a nondegenerate {\it associative form} $\lambda$, that is,
a symmetric bilinear form $\lambda:L\times L\to \F$
satisfying $\lambda([\xi,\eta],\theta)= \lambda(\xi,[\eta,\theta])$ for all
$\xi,\eta,\theta\in L$.
(Note that the latter condition together with anticommutativity of the Lie bracket easily implies that
$\lambda([\xi,\eta],\theta)=
\lambda(\theta,[\xi,\eta])$, hence the symmetry of $\lambda$ is automatic if $L$ is perfect.
In particular, in view of the interpretation of associativity which we are about to give,
there are non nonzero
$L$-module homomorphisms $L\wedge L\to\F$ if $L$ is perfect, in odd characteristic.)
Since the associativity condition can be written in the equivalent form
$\lambda([\eta,\xi],\theta)+\lambda(\xi,[\eta,\theta])=0$,
it simply means that the corresponding
linear map $\lambda:L\otimes L\to\F$ is a homomorphism of $L$-modules
into the trivial module.
Consequently, the adjoint module of $L$ is self-dual;
this condition is, in fact, equivalent with the existence
of a nondegenerate bilinear form on $L$ satisfying associativity but not necessarily symmetric.

By composition with the inverse of the $L$-module isomorphism $L\to L^{\ast}$
given by $\xi\mapsto \lambda(\xi,\cdot)$, the monomorphism
$H^2(L,\F)\to H^1(L,L^{\ast})$ turns into a monomorphism
$H^2(L,\F)\to H^1(L,L)=\Der(L)/\ad(L)$,
where $\ad(L)$ is the space of inner derivations of $L$.
Its image is $\SkDer(L)/ \ad(L)$,
where $\SkDer(L)$ denotes the space of
all derivations $D:L\to L$ which are {\it skew} with respect to the
associative form $\lambda$, that is, which satisfy
$\lambda(D(\xi),\eta)=-\lambda(D(\eta),\xi)$ for all $\xi,\eta\in L$
(see the Remark after Proposition~1.3 of~\cite{Far:central}).
Writing this condition in the equivalent form
$\lambda(D(\xi),\eta)+\lambda(\xi,D(\eta))=0$
shows that a derivation $D$ of $L$ is skew
exactly if $D$ annihilates the form $\lambda$ viewed as an element of $(L\otimes L)^\ast$,
the dual of the tensor square of the adjoint module of $L$;
it follows, in particular, that $\SkDer(L)$
is a $p$-subalgebra of $\Der(L)$ containing all inner derivations of $L$.

The isomorphism $H^2(L,\F)\to\SkDer(L)/ \ad(L)$
is actually induced by an isomorphism
$Z^2(L,\F)\to\SkDer(L)$,
which we describe here for convenience.
Because of the nondegeneracy of $\lambda$, for each cocycle $\varphi\in Z^2(L,\F)$
there is a unique derivation $D_{\varphi}:L\to L$, necessarily skew, such that
\begin{equation*}
\lambda(D_{\varphi}(\xi),\eta)=\varphi(\xi,\eta)
\qquad\text{for all $\xi,\eta\in L$.}
\end{equation*}
Conversely, the $2$-cocycle associated with the skew derivation
$D:L\to L$ is given by
\begin{equation*}
\varphi(D)(\xi,\eta)=
\lambda(D(\xi),\eta)
\qquad\text{for all $\xi,\eta\in L$.}
\end{equation*}

Now we apply these well-known facts to the Hamiltonian algebras under consideration
(the case of $H(2:\n;\omega_0)$ being already dealt with in~\cite{Far:central}).
It is known from \cite[Theorem~4.4]{Far:associative} or
\cite[Chapter~4, Theorem~6.5]{SF} that the graded algebra
$H(2:\n;\omega_0)$ has a non-degenerate associative form
$\lambda$, which in our notation becomes
\begin{equation*}
\lambda(
x^{(i)}y^{(j)},
x^{(k)}y^{(l)})=
(-1)^{i+j}\,
\delta(i+k,p^{n_1}-1)\,
\delta(j+l,p^{n_2}-1).
\end{equation*}
Thus, the dual basis of $\{x^{(i)}y^{(j)}\}$ with respect to the
nondegenerate form $\lambda$ is given by
$(x^{(i)}y^{(j)})^{\ast}=
(-1)^{i+j} x^{(p^{n_1}-1-i)}y^{(p^{n_2}-1-j)}$.

We assume now that $p>3$.
The derivation algebras of the simple Lie algebras of Cartan type
are known and are summarized in~\cite[pp.~903--905]{BKK}.
(Alternatively, the derivation algebras of Block algebras,
which include the Hamiltonian algebras under consideration here,
were already computed in~\cite[Theorem~14]{Bl}, again for $p>3$.)
In particular, it is known that
$\Der H(2:\n;\omega_i)=
\overline{C\tilde H(2:\n;\omega_i)}$,
the $p$-closure of $C\tilde H(2:\n;\omega_i)$ in
$\Der \F[2:\n]$,
for $i=0,2$.
More explicitly, a basis for the
space of outer derivations $L=H(2:\n;\omega_0)$
(or, more precisely, a set of
representatives for a basis of
$H^1(L,L)=\operatorname{Der}(L)/\ad(L)$)
consisting of homogeneous derivations with respect to the $\Z^2$-grading
is as follows.
\begin{enumerate}
\item $(\ad x)^{p^r}$, of degree $(0,-p^r)$, for $0<r<n_2$;
\item $(\ad y)^{p^s}$, of degree $(-p^s,0)$, for $0<s<n_1$;
\item $\ad (x^{(p^{n_1})})$, of degree $(p^{n_1}-1,-1)$;
\item $\ad (y^{(p^{n_2})})$, of degree $(-1,p^{n_2}-1)$;
\item $\ad(\bar x\bar y)= [\ad (x^{(p^{n_1})}),\ad (y^{(p^{n_2})})]$,
  of degree $(p^{n_1}-2,p^{n_2}-2)$;
\item the {\it degree derivation} $\ad h$,
  which has degree $(0,0)$ and acts as\\
  $(\ad h)(x^{(i)}y^{(j)})= (i+j-2)x^{(i)}y^{(j)}$.
\end{enumerate}
Note that the derivations under (1) and (2) are powers
of inner derivations of $H(2:\n;\omega_0)$;
together with the inner derivations, they span its $p$-closure $\overline{H(2:\n;\omega_0)}$.
We have denoted the derivations under (3), (4) and (5) as restrictions
of inner derivations of $\tilde H(2:\n;\omega_0)$.
Finally, the degree derivation is the restriction
of the inner derivation $\ad h$ of $W(2:\n)$, where the element
$h=x\,\partial/\partial x+y\,\partial/\partial y$
has no analogue in the Poisson bracket notation which we have adopted.

All derivations listed above except the degree derivation
(because $p>3$, but see Remark~\ref{rem:char3} for $p=3$) are skew
with respect to $\lambda$
(by direct verification, or from~\cite[Proposition~2.2]{Far:central}), and hence
$\dim(H^2(L,\F))=\dim(\Der(L)/\ad(L))-1=n_1+n_2+1$
for $L$ and $p>3$
(cf.~\cite[Theorem.~2.4]{Far:central} and~\cite{Dzhu:central}).
A set of $n_1+n_2+1$ cocycles of $L=H(2:\n;\omega_0)$
which form a basis of $H^2(L,\F)$ can be obtained from the skew derivations described above
according to the procedure described earlier.
Note that since $\lambda:L\otimes L\to\F$ is a graded map
of degree $(-p^{n_1}+3,-p^{n_2}+3)$ (where the trivial module $\F$ is assigned degree zero),
the cocycles thus obtained are homogeneous with respect to the grading
of $Z^2(L,\F)$ inherited by the $\Z^2$-grading of $L$.
A conclusion which is more relevant for us is that the universal central extension $M$ of $L$
(see Remark~\ref{rem:AFS-cocycles})
inherits a $\Z^2$-grading from $L$, and the central elements
{\em corresponding}
(see Remark~\ref{rem:graded-cocycles})
to the cocycles obtained are homogeneous.
We record the degrees explicitly:
if $D$ is a homogeneous skew derivation of degree
$(i,j)$, then the central elements of the universal central extension $M$ of $L$
corresponding to the cocycle $\varphi(D)$ acquires degree $(p^{n_1}-3-i,p^{n_2}-3-j)$.

\begin{rem}\label{rem:graded-cocycles}
Strictly speaking, central elements of $M$ correspond naturally to elements of the second homology group $H_2(L,\F)$,
and not of its dual $H^2(L,\F)$.
There is, however, a natural correspondence between the homogeneous components in the $A$-gradings of
$H_2(L,\F)$ and its dual $H^2(L,\F)$, reversing the sign of the degrees.
The fact that all these components are one-dimensional in the present case
(and in the case of $L=H(2:\n;\omega_2)$ below) justifies our abuse of language.
\end{rem}

Now we turn our attention to the filtered algebra
$H(2:\n;\omega_2)$.
This algebra has a non-degenerate associative form $\lambda$
(see~\cite[Theorem~7]{Bl}), defined by the
same formula given above for $H(2:\n;\omega_0)$ with, in addition,
$\lambda(\bar x\bar y,\bar x\bar y)=1$ and $\lambda(\cdot,\cdot)=0$
in all remaining cases.
Again from \cite{BKK} (or the original source \cite[Theorem~3.2]{Kuz:truncated}), all outer derivations
of $H(2:\n;\omega_2)$ are lifted from part of those of its
associated graded algebra, which is $H(2:\n;\omega_0)\oplus\langle\bar x\bar y\rangle$.
Specifically, a basis for the space of outer derivations of
$H(2:\n;\omega_2)$ is given by
\begin{equation*}
(\ad x)^{p^r}
\text{ for }
0<r\leq n_2,
\text{\qquad and \qquad}
(\ad y)^{p^s}
\text{ for }
0<s\leq n_1.
\end{equation*}
(Note that
$(\ad x)^{p^{n_2}}$ and
$(\ad y)^{p^{n_1}}$
induce the derivations
$\ad (x^{(p^{n_1})})$ and
$\ad (y^{(p^{n_2})})$ on the associated graded algebra.)
In particular, all derivations of $L=H(2:\n;\omega_2)$
belong to its $p$-closure in $\Der(L)$ and, consequently, they are all skew.

So far our assumption that $p>3$ was in force.
However, the $n_1+n_2$ derivations of $L=H(2:\n;\omega_2)$ which we have described
clearly remain linearly independent in $\Der(L)/\ad(L)$ also for smaller characteristics,
provided in characteristic two we replace the degree derivation
(which coincides with the inner derivation
$\ad(x\,\partial/\partial x-y\,\partial/\partial y)=
\ad(\mathcal{D}(xy))$ in that case)
with $\ad(x\,\partial/\partial x)$.
In Section~\ref{sec:AF} we will identify $H(2:\n;\omega_2)$ with an algebra of Albert and Frank
(a special type of Block algebra, see Section~\ref{sec:Hamiltonian}),
and we will prove (extending results obtained in~\cite{Bl} for characteristic $p>3$)
that its space of outer derivations has dimension $n_1+n_2$,
in every positive characteristic.
Consequently, the derivations of $H(2:\n;\omega_2)$ described above,
but with $\ad(x\,\partial/\partial x)$ replacing the degree derivation
(a change which is only relevant in characteristic two),
form a basis for its space of outer derivations regardless of the characteristic.

Setting $\varphi_r=\varphi((\ad x)^{p^r})$ and
$\psi_s=\varphi((\ad y)^{p^s})$, the following result follows by direct computation.

\begin{theorem}\label{thm:cocycles}
A basis for the second cohomology group $H^2(L,\F)$ of the
Hamiltonian Lie algebra $L=H(2:\n;\omega_2)$ over a field $\F$ of odd characteristic
is given by the classes of the cocycles $\varphi_r$ and $\psi_s$, for
$0<r\leq n_2$ and $0<s\leq n_1$, as defined by the following formulas:
\begin{align*}
\varphi_r(
x^{(i)}y^{(j)},
x^{(k)}y^{(l)})
&=
\begin{cases}
  (-1)^{i+j} & \text{if
    $(i+k,j+l)\equiv (0,p^r)$}\\
  & \text{ \qquad$\pmod{(p^{n_1}-1,p^{n_2}-1)}$}
  \\
  0 & \text{otherwise;}
\end{cases}
\\
\psi_s(
x^{(i)}y^{(j)},
x^{(k)}y^{(l)})
&=
\begin{cases}
  (-1)^{i+j+1} & \text{if
    $(i+k,j+l)\equiv (p^s,0)$}\\
  & \text{ \qquad$\pmod{(p^{n_1}-1,p^{n_2}-1)}$}
  \\
  0 & \text{otherwise.}
\end{cases}
\end{align*}
\end{theorem}

The $A$-degrees of the cocycles $\varphi_r$ and $\psi_s$ are
$(2,-p^r+2)$ and $(-p^s+2,2)$, respectively.
The central elements of the universal central extension $M$ of $L=H(2:\n;\omega_2)$
(see Remark~\ref{rem:AFS-cocycles})
corresponding to them acquire $A$-degree
$(-2,p^r-2)$ and $(p^s-2,-2)$, respectively.

\begin{rem}\label{rem:char3}
It is easy to verify that the degree derivation of $H(2:\n;\omega_0)$ acts on $\lambda$,
considered as an element of $(L\otimes L)^{\ast}$,
as multiplication by $6$.
In particular, when the characteristic is three (or two, but see Remark~\ref{rem:char2-cocycles} concerning this case)
the degree derivation is skew
(cf.~\cite[Proposition~2.2]{Far:central}),
and the second cohomology group of $H(2:\n;\omega_0)$ becomes larger.
Also, the derivation algebra of $L$ can be larger (see~\cite[p.~197]{SF}).
In particular, the graded Hamiltonian algebra
$L=H(2:(1,1);\omega_0)$ in characteristic three is a classical Lie algebra of type $A_2$,
namely, $L$ is isomorphic to the quotient of $\Sl_3$ modulo its one-dimensional center,
see~\cite[Lemma~6.4]{Skr:low}.
According to~\cite[Corollary~3]{Bl} and the identification of Hamiltonian algebras with Block algebras,
this is the only instance for $p>2$ where an algebra $H(2:\n;\omega_0)$
or $H(2:\n;\omega_2)$ is isomorphic with a classical algebra.
It is well known that $\Der(L)$ is a fourteen-dimensional classical algebra of type $G_2$
(cf.~\cite[p.~678]{Skr:low}), hence
$\dim(\Der(L)/\ad(L))=7$.
The additional three derivations with respect to those described earlier can be obtained from
$\ad x^{(3)}$, $\ad y^{(3)}$ and $\ad\bar x\bar y$ by conjugation under the automorphism
$x^{(i)}y^{(j)}\mapsto
(-1)^{ij}\, x^{(2-i)}y^{(2-j)}$
of $L$.
Since all derivations are skew in this case, we conclude that
$H^2(L,\F)$ has dimension $7$
(but see \cite[p.~38]{vdK} for another proof of this fact).

It is easily checked that exactly one of
those three `exceptional derivations' of
$H(2:(1,1);\omega_0)$ extends to a derivation of $H(2:(1,n_2);\omega_0)$ for $n_2>1$.
These additional derivations account for all the exceptions in characteristic three with respect to the
description of outer derivations of $H(2:\n;\omega_0)$ in higher characteristic
given earlier in this section.
This follows from~\cite[Proposition~4.3]{Skr:low} and a general fact about derivations of nonnegative degree
(in the standard grading) of graded Lie algebras of Cartan type, see e.g.~\cite[Chapter~4, Proposition~8.3]{SF}.
The derivations of $L=H(2:\n;\omega_0)$ in characteristic three can be summarized as follows:
$\dim(\Der(L)/\ad(L))$ equals $n_1+n_2+2$ if $n_1,n_2>1$
(like in higher characteristic),
it equals $n_1+n_2+3$ if $n_1=1<n_2$,
and it equals $7$ if $n_1=n_2=1$.
Since all derivations of $L$ are skew in characteristic three, we have
$\dim(H^2(L,\F))=\dim(\Der(L)/\ad(L))$.
\end{rem}

\begin{rem}\label{rem:char2-derivations}
Again according to~\cite[Proposition~4.3]{Skr:low} and~\cite[Chapter~4, Proposition~8.3]{SF},
the derivations of a simple algebra $H(2:\n;\omega_0)$ in characteristic two (hence with $n_1,n_2>1$)
allow the same description
as in characteristic greater than three, by the list given earlier in this section,
except that the degree derivation in item (6) (which is inner in characteristic two, as it coincides with $\ad xy$)
should be replaced with the derivation acting as $D(x^{(i)}y^{(j)})=(i-1)\,x^{(i)}y^{(j)}$.
In particular, $\Der(L)/\ad(L)$ has dimension $n_1+n_2+2$.

The algebra $H(2:(1,n_2);\omega_0)$ in characteristic two has more derivations than usual.
Since it is a semidirect product of a simple Zassenhaus algebra by its adjoint module, its derivations
can be easily calculated from those of the Zassenhaus algebra (see Remark~\ref{rem:Zas-der-cocycles}).
In fact, if $L$ is a semidirect product of a simple algebra $S$ by its adjoint module, then
$\dim(\Der(L))=2\,\dim(\Der(S))+2$ in characteristic two, and
$\dim(\Der(L))=2\,\dim(\Der(S))+1$ otherwise.
Consequently, $\Der(L)/\ad(L)$ has dimension $2\,n_2+2$ for $L=H(2:(1,n_2);\omega_0)$.
Alternatively, since $H(2:(1,n_2);\omega_0)$ is the tensor product of a simple Zassenhaus algebra
with a ring of divided powers $\F[z:1]$, the conclusion follows from~\cite[Theorem~7.1]{Bl:differentiably-simple}.
\end{rem}

\begin{rem}\label{rem:char2-cocycles}
In characteristic two the argument which relates derivations of $L$ into $L^\ast$ and
the second cohomology group of $L$ needs to be modified as follows.

The image of the map $Z^2(L,\F)\to Z^1(L,L^\ast)$ consists of all derivations which are
{\em alternating}, in the sense that $D(\xi)(\xi)=0$ for all $\xi\in L$.
This condition is equivalent to being skew in odd characteristic, but is stronger in characteristic two.
(It is convenient to reserve the term {\em skew} for the weaker condition, as it
applies in slightly greater generality;
see~\cite{Far:central}, where Lemma~1.1 remains valid for skew derivations, but not for alternating derivations,
in characteristic two.)
In presence of a nondegenerate associative form $\lambda$,
a derivation $D$ of $L$ will be called {\em alternating} with respect to $\lambda$ if
$\lambda(D(\xi),\xi)=0$ for all $\xi\in L$.
Since inner derivations of $L$ are alternating, there is an isomorphism of $H^2(L,\F)$ with
the quotient of the space of alternating derivations
by the space of inner derivations of $L$.
Note that while the space of skew derivations is a $p$-subalgebra of $\Der(L)$
(in every characteristic), the space of alternating derivations is a Lie subalgebra
but need not be a $p$-subalgebra in characteristic two (as shown by the examples below).

In order to verify that a skew derivation is alternating it suffices to check that
$D(\xi)(\xi)=0$ for all elements $\xi$ of some basis of $L$.
Also, if $L$ is graded and its associative form is homogeneous with respect to the grading
(as is the case for the $\Z^2$-grading of $H(2:\n;\omega_0)$ and the $A$-grading of $H(2:\n;\omega_2)$)
it is enough to check derivations which are homogeneous with respect to the grading,
because the alternating derivations
form a graded subalgebra (as well as the skew derivations).

For $L=H(2:\n;\omega_0)$ with $n_1,n_2>1$, all derivations described earlier in this section
(taking Remark~\ref{rem:char2-derivations} into account) are alternating, except $\ad x^{(p^{n_1})}$
and $\ad y^{(p^{n_2})}$, which are only skew.
Therefore, $H^2(L,\F)$ has dimension $n_1+n_2$.
An examination of $L=H(2:(1,n_2);\omega_0)$ shows that the alternating derivations
are exactly those which normalize the Zassenhaus subalgebra $\langle xy^{(j)}:j=0,\ldots,2^{n_2}-2\rangle$
(together with the inner derivations).
These correspond to the derivations described under items (1), (5) and (6)
in the list given earlier in this section,
and we conclude that $H^2(L,\F)$ has dimension $n_1+n_2$ in this case, too.

In the case of $L=H(2:\n;\omega_2)$
all alternating derivations are inner.
As a consequence, in characteristic two we have $H^2(L,\F)=0$.
\end{rem}

\begin{rem}\label{rem:Zas-der-cocycles}
Recalling from Section~\ref{sec:Hamiltonian} that the simple Zassenhaus algebra $L=W(1:n)^{(1)}$
in characteristic two is isomorphic with $L=H(2:(1,n-1);\omega_2)$,
the previous remark shows that its second cohomology group $H^2(L,\F)$ vanishes.
By contrast, the second cohomology group of the Zassenhaus algebra $W(1:n)$ in odd characteristic has
dimension one if $p>3$, and dimension $n-1$ if $p=3$
(as a special case of~\cite[Theorem.~3.2]{Far:dual} or~\cite{Dzhu:central}).
The second cohomology group of the simple Zassenhaus algebra was also computed
in~\cite[Theorem~2]{Dzhu:central-Zassenhaus};
however, note that the 
central extensions of $W(1:n)^{(1)}$ in characteristic two
which are exhibited there are not Lie algebras in the common sense, because their multiplication
is (skew-)symmetric but not alternating.
For the sake of completeness we mention that the algebra of outer derivations of the simple
Zassenhaus algebra $W(1:n)^{(1)}$ has dimension $n-1$ if $p$ is odd and $n$ if $p=2$.
This is well known, but the case where $p=2$ is also a consequence of Theorem~\ref{thm:AF-derivations}.
\end{rem}

\begin{rem}\label{rem:G2}
The last sentence in~\cite{Bl} claims that there is no isomorphism between a Block algebra
and an exceptional (classical) simple algebra except when $p=2$ and the algebra has dimension 14,
since otherwise their dimensions are distinct.
This may leave some doubt on whether in characteristic two $H(2:(2,2);\omega_0)$ might be isomorphic with a simple
algebra of type $G_2$ (which, in turn, is isomorphic with the quotient of $\Sl_4$ by its one-dimensional center).
However, according to~\cite{vdK} the second cohomology group of the latter has dimension $7$,
and this fact together with Remark~\ref{rem:char2-cocycles} excludes the possibility of an isomorphism.
\end{rem}

\section{Some cyclic gradings of $H(2:\n;\omega_2)$}
\label{sec:gradings}

The $A$-grading of $H(2:\n;\omega_2)$ defined in Section~\ref{sec:Hamiltonian}
leads in a natural
way to several gradings (here called {\em specializations}) over cyclic quotients $\bar A$ of
$A=\Z^2/\langle(p^{n_1}-1,p^{n_2}-1)\rangle$.
More precisely, for any pair of integers $(R,S)$ and any divisor $N$ of $R(p^{n_1}-1)+S(p^{n_2}-1)$
we have a group homomorphism
$\mu:  A \to \bar A= \Z/N\Z$ given by
$\mu(i,j)=Ri+Sj+N\Z$.
Correspondingly, we obtain an $\bar A$ grading $L=\bigoplus_{k\in \bar A} L_k$ by setting
$L_k=\sum_{\mu(i,j)=k} L_{i,j}$.
In what follows we set
$N=|R(p^{n_1}-1)+S(p^{n_2}-1)|$,
since the remaining cases can be obtained from these gradings through a further specialization.
Also, it is no loss to assume that the homomorphism is surjective,
which amounts to choosing $R$ and $S$ relatively prime (because of our choice of $N$).
To help visualizing the grading thus obtained, it may be convenient to arrange the monomials
in $H(2:\n;\omega_2)$ in a $(p^{n_1}\times p^{n_2})$-array according to the degrees of $x$ and $y$,
and think of the specialization process as {\em slicing} the $A$-grading
according to some specified direction.

The simplest specialization where $(R,S)=(0,-1)$ will be useful in Section~\ref{sec:thin}.
This is a $\Z/N\Z$-grading with $N=p^{n_2}-1$, every
component has dimension $p^{n_1}$ and is spanned by all monomials
where $y$ has a given degree $1\leq j\leq p^{n_2}-2$, except the
component of degree $1$, which has dimension $2p^{n_1}-1$ and is
spanned by all monomials where $y$ has degree $0$ or $p^{n_2}-1$.
The component of degree $0$ is isomorphic with a Zassenhaus algebra $W(1:n_1)$.

The following two specializations of the $A$-grading of $H(2:\n;\omega_2)$ are more interesting.

\subsection{A grading related to graded Lie algebras of maximal class}\label{subsec:maximal}
Let $(R,S)=(-p^{n_2},-1)$.  Then $N=p^n-1$, and all
components $L_k$ are one-dimensional.  In fact,
$L_{ip^{n_2}+j}=\langle x^{(p^{n_1}-i)} y^{(p^{n_2}-j)} \rangle$, for
$0<i\leq p^{n_1}$ and $0<j\leq p^{n_2}$ with
$(i,j)\not=(p^{n_1},p^{n_2})$.
Furthermore, $L$ has an outer derivation
\begin{equation}\label{eq:nonsing-der}
D=\bar y\frac{\partial}{\partial x}+(1+e)\frac{\partial}{\partial y},
\end{equation}
which is non-singular and homogeneous of degree one (with respect to the grading under consideration).
Hence $D$ permutes the components of the grading cyclically, namely $DL_k=L_{k+1}$ for all $i$.
(Recall from Section~\ref{sec:Hamiltonian} that we write $D$ for $\ad D$,
since the latter acts the same way as a Lie algebra derivation of
$P(2:\n;\omega_2)$ as $D$ acts as a derivation of the associative structure.)
In fact,
we have
\begin{equation*}
Dx^{(i)}
y^{(j)}=
\bar y x^{(i-1)}
y^{(j)}
+
(1+e)
x^{(i)}
y^{(j-1)}
=
\begin{cases}
  x^{(i-1)} \bar y &\text{if $j=0$,}
  \\
  (1+e) x^{(i)} y^{(j-1)} &\text{if $j>0$.}
\end{cases}
\end{equation*}
In particular, the derivation $D$ is periodic of period $p^n-1$.
Note that $D$ is the derivation denoted by $D_2$ in \cite[p.~911]{BKK},
viewing the Hamiltonian algebra $H(2:\n;\omega_2)$
as the special algebra $S(2:\n;\omega_2)$.

We quote from~\cite{Kos:beginnings} the following definition.

\begin{definition}\label{def:maximal}
We say that a finite-dimensional Lie algebra $L$ admits a nonsingular derivation $D$
{\em agreeing with a $\Z/N\Z$-grading} $L=\bigoplus_{k\in\Z/N\Z} L_k$
if $DL_k=L_{k+1}$ for all $k\in\Z/N\Z$.
\end{definition}

It will be convenient to allow any finite cyclic group to replace $\Z/N\Z$ in the definition,
provided we specify a distinguished generator of it (to play the role of $1$).
The situation described in Definition~\ref{def:maximal} where all components $L_k$
have dimension one
played a crucial role in~\cite{ShZe:finite-coclass} and~\cite{Sha:coclass}.

Suppose $L$ is a finite-dimensional Lie algebra possessing a nonsingular derivation
which agrees with a $\Z/N\Z$-grading with one-dimensional components.
We build the corresponding {\em twisted loop algebra}
$\bigoplus_{k\in\Z} L_{\bar k}\otimes t^k$ inside $L\otimes_{\F}\F[t,t^{-1}]$,
where $\bar k$ denotes the residue class of $k$ modulo $N$.
The Lie algebra spanned by its positive part
$\bigoplus_{k>0}L_{\bar k}\otimes t^k$
together with its derivation $D\otimes t$
is a graded Lie algebra of maximal class
in the sense of~\cite{CMN}.
With a harmless abuse of language we will call the latter
{\em the loop algebra} of $L$.

In particular, the loop algebra of
$H(2:\n;\omega_2)$ with respect to the derivation $D$
and the grading which we have just constructed is a graded Lie algebra of
maximal class, and precisely one of those which we have named after
Albert-Frank-Shalev in \cite{CMN}.
We will come back to this grading in Section~\ref{sec:AF}.

\begin{rem}\label{rem:AFS-cocycles}
We comment briefly on the relevance of the second cohomology group
of $H(2:\n;\omega_2)$, which we have discussed in Section~\ref{sec:cohomology}, to presentations
of the algebras of Albert-Frank-Shalev $AFS(a,b,n,p)$ (see Section~\ref{sec:AF} for their definition).
Although these algebras are not finitely presented,
it is proved in~\cite{Carrara} that they
are quotients of certain finitely presented Lie algebras modulo their second centers.
Knowledge of the second cohomology group of $H(2:\n;\omega_2)$ sheds light
on these particular extensions of the algebras of Albert-Frank-Shalev, as we illustrate below.
The second cohomology group of $H(2:\n;\omega_0)$ plays a similar role in~\cite{Avi} and~\cite[Section~4]{CaMa:Nottingham}.

Recall that every perfect Lie algebra $L$ has a universal central extension
$0\to Z\to M\to L\to 0$
(see~\cite[Section~1]{vdK}).  
In particular, $M/Z\cong L$, so we may view $L$ as a quotient of $M$,
and $Z\cong H_2(L,\F)\cong H^2(L,\F)^\ast$ as vector spaces.
According to~\cite[Theorem~2.2]{BenMoo}, every derivation of $L$
lifts to a derivation of $M$;
if $L$ is centerless, the lift is unique, therefore $\Der(M)\cong\Der(L)$,
and $M$ is a $\Der(L)$-module in a natural way.
In the case of $L=H(2:\n;\omega_2)$, we claim that $\Der(L)$ acts trivially on $Z$,
which coincides with the center $\zeta(M)$ of $M$ here.
In fact, for a perfect and centerless Lie algebra $L$,
the dual module $\zeta(M)^\ast$ and $H^2(L,\F)$ are easily seen to be isomorphic
not only as vector spaces, but as $\Der(L)$-modules.
Now assume, in addition, that the characteristic is odd and that $L$ has a nonsingular associative form $\lambda$.
Then the isomorphism of
$H^2(L,\F)$ with $\SkDer(L)/\ad(L)$ described in Section~\ref{sec:cohomology} is also an isomorphism of $\Der(L)$-modules
(with respect to the adjoint action of $\Der(L)$ on itself).
Since $\Der(L)/\ad(L)$ is abelian for $L=H(2:\n;\omega_2)$, our claim follows.
Therefore, the center of the extension of $M$ by $\langle D\rangle$,
where $D$ is the nonsingular derivation~\eqref{eq:nonsing-der}, coincides with $\zeta(M)$ and, in particular, is nonzero
according to Theorem~\ref{thm:cocycles}, in odd characteristic.
It follows that the loop algebra $\tilde M$ of $M$ with respect to $D$ has an infinite-dimensional center.
The quotient of $\tilde M$ by its center is isomorphic with $\tilde L$.
A standard result of B.~H.~Neumann recalled in~\cite{CaMa:thin} as Theorem~6
implies that the quotient of a finitely generated Lie algebra modulo an infinite-dimensional central ideal
cannot be finitely presented; in particular,
$\tilde L\cong\tilde M/\zeta(\tilde M)$ is not finitely presented.

The main result of~\cite{Carrara} shows that a suitable central extension of $\tilde M$ is finitely presented.
(The need to take a further central extension to obtain a finitely presented algebra is due to the fact
that the second cohomology group of a loop algebra, besides depending on the second cohomology
group of the underlying finite-dimensional algebra, includes a component arising from associative forms
of the latter and the cyclic homology of the polynomial ring $\F[t]$ which we are tensoring with.
We will not discuss this point further here, but see~\cite{Zus:current-central}.)
According to our identification in Theorem~\ref{thm:AF=Hamiltonian} of the algebras of Albert and Frank
with Hamiltonian algebras $H(2:\n;\omega_2)$, the second cohomology group of the latter
discussed in Theorem~\ref{thm:cocycles} can be recognized in (part of) the central elements
of the finitely presented central extensions of the Albert-Frank-Shalev algebras
considered in~\cite{Carrara}.
More precisely, the $\bar A$-grading of $L=H(2:\n;\omega_2)$ (like any other grading)
extends uniquely to a grading of its universal central extension $M$.
The central elements of the latter corresponding (recall~Remark~\ref{rem:graded-cocycles}) to the cocycles
$\varphi_r$ and $\psi_s$ of Theorem~\ref{thm:cocycles} (in odd characteristic) occur in degrees
$2q-p^r+2$ and $2q-qp^s+2$, respectively,
and give rise to central elements of $\tilde M$, the loop algebra of $M$ with respect to $D$,
in all degrees congruent to these modulo $\dim(L)=p^{n_1+n_2}-1$.
These central elements can be recognized in the list given in~\cite[p.~399-400]{Carrara}, in the special case
of the algebra $AFS(a,n,n,p)$.
(Additional complications arise in the case of characteristic two, which we have neglected here for simplicity.)
\end{rem}

By symmetry, the case $(R,S)=(-1,-p^{n_1})$ is completely analogous.
In this case
\begin{equation*}
D^{p^{n_2}}=(1+e)\frac{\partial}{\partial x}+\bar x\frac{\partial}{\partial y}
\end{equation*}
is a non-singular derivation which permutes the components cyclically.
This latter grading is just one instance of a whole set of $\bar A$-gradings which can be obtained
from the former $\bar A$-grading by an application of an automorphism of the grading group $\bar A$.
More precisely,
$L$ is also graded by the subspaces $\tilde L_i = L_{ki}$, where $k$ is any integer with $(k,p^n-1)=1$.
In general, a derivation $\tilde D$
such that $\tilde D\tilde L_i= \tilde L_{i+1}$
need not exist.
However, it certainly does if  $k$ is a power of $p$, because then
$D^{p^t}\tilde L_i=\tilde L_{i+1}$,
if $\tilde L_i = L_{p^ti}$.
This way of obtaining new gradings is related to the process of {\em deflation}
for graded Lie algebras of maximal class introduced in~\cite{CMN}.
In fact, the loop algebra of $L$ with respect to its grading
given by the subspaces $\tilde L_i=L_{pi}$ and its derivation $D^p$
is the deflation of
the loop algebra of $L$ with respect to its grading
given by the subspaces $L_i$ and its derivation $D$.

\subsection{A thin grading}\label{subsec:thin}
Let $(R,S)=(-p^{n_2}+1,-1)$.
This is a $\Z/N\Z$-grading with $N=p^{n_1}(p^{n_2}-1)$,
and the components have dimension one or two.
In the present grading, the two-dimensional components
are those of degree $i(p^{n_2}-1)+1$ for $1<i\leq p^{n_1}$.
In particular, $L_1=\langle x, \bar y\rangle$.
This grading fits the following definition.

\begin{definition}\label{def:thin-fd}
A grading $L=\bigoplus_{k\in\Z/N\Z} L_k$ of a (finite-dimensional) Lie algebra $L$
over a field $\F$ is called {\em thin} if $L_1$ is two-dimensional,
and the following {\it covering property\/} holds
\[
 \text{for all $k\in\Z/N\Z$, and all $u \in L_{k}$, $u \ne 0$,
 we have $L_{k+1} = [u, L_{1}]$}.
\]
\end{definition}

Again, the definition is motivated by an analogous one for positively graded,
possibly infinite-dimensional Lie algebras.
In fact, given a thin $\Z/N\Z$-grading of a finite-dimensional Lie algebra $L$,
the (positive part of the twisted) loop algebra
$\bigoplus_{k>0}L_{\bar k}\otimes\F t^k$
is a thin Lie algebra
in the sense of the following definition (see~\cite{CaMa:thin} for background).

\begin{definition}\label{def:thin}
A graded Lie algebra $L=\bigoplus_{k=1}^{\infty} L_k$
is called {\em thin} if $L_1$ is two-dimensional,
and the following {\it covering property\/} holds
\[
 \text{for all $k\ge 1$, and all $u \in L_{k}$, $u \ne 0$,
 we have $L_{k+1} = [u, L_{1}]$}.
\]
\end{definition}

With both definitions, it follows from the covering property that
$L$ is generated by $L_1$, and that
$\dim(L_{k}) \le 2$ for all $k$. 
We call a homogeneous component $L_{k}$ of dimension $2$ a {\it diamond}.
The diamonds will be numbered in the natural order of occurrence
(cyclic starting from $L_1$ in case of Definition~\ref{def:thin-fd}).
Therefore, $L_{1}$ is the first diamond.
If there are no other diamonds in case of Definition~\ref{def:thin},
then $L$ is an algebra of maximal class (see \cite{CMN}).
We refer to
the finite sequence of one-dimensional homogeneous components between two
consecutive diamonds as a {\it chain}.

We must point out that there are also instances where a thin Lie algebra
in the sense of Definition~\ref{def:thin} is constructed as a loop algebra
from a suitable grading of a finite dimensional simple (Hamiltonian) Lie algebra which
does not quite fit Definition~\ref{def:thin-fd}, but requires the intervention of a nonsingular
outer derivation, very much like in the construction of Lie algebras of maximal class
described in the previous subsection (see~\cite{Avi, Avi:thesis}).
These algebras are not needed in this paper, however.

The fact that an algebra with a thin grading in the sense of Definition~\ref{def:thin-fd}
gives rise to a thin Lie algebra in the sense of Definition~\ref{def:thin}
via the loop algebra construction
allows one to apply to the former setting arguments and results originally formulated for the latter.
For example, results from~\cite{CMNS}, extended in~\cite{AviJur},
imply that in a finite-dimensional Lie algebra $L$ over a field
of arbitrary characteristic
with a thin grading, the second diamond
can only occur in degree $3$, $5$, $q$, or $2q-1$,
where $q$ is a power of the characteristic.
Some care is needed in carrying definitions over from the infinite-dimensional setting to the present one,
where the degree of a homogeneous element is an integer defined only modulo the order $N$ of the grading group:
when speaking of {\em the degree where the second diamond occurs,}
we actually refer to the smallest integer greater than one in which degree a diamond occurs.

Suppose now that $L$ is a finite-dimensional Lie algebra with a thin grading,
and suppose that the second diamond occurs in degree
$2q-1$, where $q$ is a power of the characteristic.
According to~\cite{CaJu:quotients} we have
$$
 C_{L_{1}} (L_{2}) = C_{L_{1}} (L_{3}) =
 \dots = C_{L_{1}} (L_{2q-3}) = \langle Y\rangle,
$$
provided the characteristic is odd;
this fails in characteristic two, as is shown in
\cite{Ju:quotients,JuYo:quotients}.
We let
$X
\in L_{1}
\setminus
\langle Y\rangle$, so that $X$ and $Y$ generate $L$.
(Note that with respect to the analogous situation in~\cite{CaMa:thin}
we have switched to capital letters for the generators $X$ and $Y$,
to avoid conflict of meaning with the variables we use for divided powers.)

Suppose that $L_k$ is any diamond of $L$, and let $L_{k-1}=\langle V\rangle$.
It is not difficult to show, as in~\cite{CaMa:thin}, that
$[V,X,Y]+[V,Y,X]=[V,Y,Y]=0$, and to deduce that
$$
 [V,Y,X] = \lambda [V,X,X]
$$
for some $\lambda \in \F \cup \{ \infty \}$,
to be read as $[V,X,X]=0$ when $\lambda=\infty$.
As in the infinite-dimensional setting of~\cite{CaMa:thin}
we will say that the diamond in degree $k$ has type $\lambda$.
Note that $\lambda$ depends on the choice made for the generators $X$ and $Y$, but
the type being finite, or infinite, does not.
Strictly speaking, $\lambda=0$ cannot occur here,
because the covering property would imply that
$[V,Y]=0$ and so $L_k=\langle[VX]\rangle$, contradicting the assumption that $L_k$ is a diamond.
However, there are situations (here and in other papers, like~\cite{CaMa:thin}, \cite{Avi} and \cite{CaMa:Nottingham})
where we have found natural
and convenient to informally call {\em fake} diamonds
certain one-dimensional components $L_k$, which thus may be assigned type $0$.
This usually happens when the diamonds of some algebra with a thin grading
(or some infinite-dimensional thin algebra) occur at regular intervals
provided we include some {\em fake} diamonds.
Unfortunately, 
which cocycle of the algebra gives rise to the central element in
a fake diamond depends on the grading under consideration, and appears not to admit an intrinsic characterization,
in terms of the algebra alone.

The $\Z/N\Z$-grading of $H(2:\n;\omega_2)$ under consideration in this subsection is a thin grading,
with diamonds in all degrees congruent to $1$ modulo $q-1$, with the exception of degree $q$,
where we have set $q=p^{n_2}$.
In fact, here we may take $Y=\bar y$ and $X=x$, and the following computations
show both the validity of the covering property and the location of the diamonds:
\begin{align*}
  & \{x^{(i)} y^{(j)}, x\} = (1+e) x^{(i)}
  y^{(j-1)},\\
  & \{x^{(i)} y^{(j)}, \bar y\} = -x^{(i-1)} y^{(j)} y^{(p^{n_2}-2)} =
  \begin{cases}
    -x^{(i-1)} y^{(p^{n_2}-2)} &\text{if $j=0$,}
    \\
    x^{(i-1)} \bar y &\text{if $j=1$,}
    \\
    0 &\text{otherwise.}
  \end{cases}
\end{align*}
Taking $V=x^{(i)}y$, for $0<i<p^{n_1}$, we see that all diamonds have type $\infty$.
Note that the elements $x^{(i)}y$ with $0\le i<p^{n_1}$, that is, the elements
{\em just above} (in the sense of {\em immediately preceding}) the diamonds,
if we include a {\em fake} diamond $L_{p^{n_2}}=\langle y\rangle$,
span a subalgebra of $H(2:\n;\omega_2)$
isomorphic with a Zassenhaus algebra $W(1:n_1)$.
This feature will reappear in the different grading which we will consider in Section~\ref{sec:thin}.

According to~\cite{CaMa:thin}, the loop algebra of $H(2:\n;\omega_2)$ with respect to this grading
is the only thin maximal subalgebra of
the algebra of maximal class $AFS(0,n_2,n,p)$, which is the loop
algebra of $H(2:\n;\omega_2)$ with respect to the grading seen
in the previous subsection. 
The grading of $H(2:\n;\omega_2)$ extends uniquely to a grading of its universal central
extension.
The central elements of the latter corresponding to the cocycles
$\varphi_r$ and $\psi_s$ of Theorem~\ref{thm:cocycles} (in odd characteristic) occur in degrees
$2q-p^r$ and $2q-p^s(q-1)$, respectively.

\section{The Hamiltonian algebra $H(2:\n;\omega_2)$ as a Block algebra}\label{sec:AF}

Benkart, Kostrikin and Kuznetsov proved in \cite[Theorem~4.9]{BKK}
(using the classification of the modular simple Lie algebras of characteristic $p>7$
completed in~\cite{Strade6})
that the only simple Lie
algebras over an algebraically closed field of characteristic $p>7$
which admit a nonsingular derivation agreeing with a $\Z/N\Z$-grading with one-dimensional components
(in the terminology introduced in Definition~\ref{def:maximal}) are of type $H(m:\n;\omega_2)$.
(In~\cite[Theorem~4.9]{BKK} $N$ had the form $p^n-1$, but this was immaterial.)
Later it was proved in~\cite{KoK3} that $m$ necessarily equals two.

Simple Lie algebras with this property were considered earlier by Shalev in \cite{Sha:max},
who noted that certain Lie algebras introduced by Albert and Frank in~\cite{AF}
can be defined over the prime field $\F_p$ and enjoy the
property.
Shalev used them to build the first examples of insoluble graded Lie algebras
of maximal class, as loop algebras of the algebras of Albert and Frank.
The algebras of Albert and Frank and their loop algebras have been further discussed in
\cite{CMN}, to which we conform our notation, as starting points for the construction
of more graded Lie algebras of maximal class.
A byproduct of the classification of graded Lie algebras of
maximal class achieved in~\cite{CN} and~\cite{Ju:maximal}
(for odd and even characteristic, respectively) is a proof,
independent of the classification of modular simple Lie algebras,
that the only finite-dimensional Lie algebras
which admit a nonsingular derivation agreeing with a $\Z/N\Z$-grading with one-dimensional components,
where $N$ is any integer prime to the characteristic,
are the algebras of Albert and Frank.
If we drop the condition on $N$ but assume the simplicity of the algebra
the result remains true in odd characteristic (extending the cited result of~\cite{BKK}),
while in characteristic two there is exactly one further class of algebras
joining the algebras of Albert and Frank
and consists of the Lie algebras constructed in~\cite{Ju:Bi-Zassenhaus} and named {\em Bi-Zassenhaus algebras}
(for which $N$ is two less than a power of two).
Note that in characteristic two the simple Zassenhaus algebra of dimension $2^n-1$ has this property,
the nonsingular derivation being given by $\ad(e_0)=\ad{(E_{-1}+E_{2^n-2})}$
in the notation introduced in Section~\ref{sec:Cartan}.
In fact, with respect to its basis $\{e_{\alpha}\}$ the simple Zassenhaus algebra
coincides with the algebra of Albert and Frank $AF(0,1,n,2)$ defined below.
(As we have mentioned in Section~\ref{sec:Hamiltonian}, it also coincides with $H(2:\n;\omega_2)$.)

The above results indirectly imply that for $p>7$ the class of algebras of Albert and Frank considered by Shalev
coincides with the class of Hamiltonian algebras $H(2:\n;\omega_2)$.
In fact, since the algebras of Albert and Frank are Block algebras this result holds for $p>3$
according to \cite[Lemma~1.8.3]{BlWil:rank-two}
together with the remarks on the characteristic which we have made at the end of Section~\ref{sec:Hamiltonian}.
In Theorem~\ref{thm:AF=Hamiltonian} we use the cyclic grading and the nonsingular derivation
of the algebras of Albert and Frank employed by Shalev
to find an explicit isomorphism of the latter with Hamiltonian algebras $H(2:\n;\omega_2)$,
in arbitrary (positive) characteristic.
We compute the derivations of the algebras of Albert and Frank in Theorem~\ref{thm:AF-derivations},
extending to arbitrary characteristic a result obtained by Block in~\cite{Bl} for $p>3$;
in view of Theorem~\ref{thm:AF=Hamiltonian},
the description of derivations of algebras $H(2:\n;\omega_2)$ quoted in Section~\ref{sec:cohomology} for $p>3$
is also extended to arbitrary characteristic.

For integers $0\le a<b< n$
(the choice made in~\cite{CMN}, or for the equivalent choice $0< a<b\le n$ made later in~\cite{CN}),
the algebra of Albert and Frank $AF(a,b,n,p)$
is the Lie algebra over $\F_{p^n}$
with basis
$\{e_{\xi}\mid \xi\in\F_{p^n}^{\ast}\}$
and multiplication given by
\begin{equation}\label{eq:AF}
[e_{\xi},e_{\eta}]=
\left(
\xi^{p^a}\eta^{p^b}-
\xi^{p^b}\eta^{p^a}
\right)
e_{\xi+\eta}
\end{equation}
(where the right-hand side is interpreted as zero when $\xi+\eta=0$,
for example by setting $e_0=0$).
Thus, an algebra of Albert and Frank is a special instance of a Block algebra (see Section~\ref{sec:Hamiltonian}).
It was shown in \cite{Sha:max} that the
algebra of Albert and Frank $S=AF(a,b,n,p)$ is defined over the prime field $\F_p$.
Moreover, it has a graded basis
$u_0,\ldots,u_{p^n-2}$ over $\F_p$, and a derivation $D$
such that $Du_i=u_{i+1}$ for all $i$, to be read modulo $p^n-1$.
Explicitly, this basis is related to the original basis
by the formulas
$u_i=\sum_{\xi\in\F_{p^n}}
\xi^{i-p^a-p^b}e_{\xi}$.
Part of the multiplication table with respect to the new basis, namely, a description of the adjoint action
of the element $u_{p^a+p^b}$, was computed in
Proposition~2.4 of~\cite{Sha:max} (but see also \cite{CMN}).
If we include the derivation $D$ into consideration
we obtain the following statement:
the split extension of $S$ by $\F_p D$ has a finite
presentation on the generators $u_0,\ldots,u_{p^n-2}$
(with subscripts viewed modulo $p^n-1$) and $D$, with
relations
\begin{equation*}
 \begin{cases}{}
   Du_i=u_{i+1} \\ 
   {}[u_{p^a+p^b}, u_{p^a}] = -u_{2p^a+p^b}, \\ 
   {}[u_{p^a+p^b}, u_{p^b}] = u_{p^a+2p^b}, \\ 
   {}[u_{p^a+p^b}, u_j] = 0 & \text{otherwise.}
 \end{cases}
\end{equation*}
This implies at once that $S$ is defined over the prime field.
Note that the above formulas differ in sign from those given in~\cite[p.~4028]{CMN},
which were incorrectly quoted from~\cite[Proposition~2.4]{Sha:max}.
That mistake amounts to using the algebra with the opposite multiplication, and
caused no serious consequence in~\cite{CMN}.
Note also that here we have written the derivation $D$ on the left, differently from~\cite{CMN}, and hence
$[u_i,D]=-[D,u_i]=-Du_i$.

In the terminology introduced in Definition~\ref{def:maximal},
the nonsingular derivation $D$ agrees with
the $\Z/N\Z$-grading $S=\bigoplus_{k\in\Z/N\Z} S_k$,
where $N=p^n-1$ and $S_k=\langle u_k\rangle$.
We will take advantage of the similar grading of $H(2:\n;\omega_2)$
defined in Subsection~\ref{subsec:maximal} to construct an isomorphism between
$S=AF(a,b,n,p)$ and $H(2:\n;\omega_2)$, where $\n=(n-b+a,b-a)$.
Since the map $e_{\xi}\mapsto e_{\xi^{p^a}}$
gives an isomorphism of $AF(a,b,n,p)$ with $AF(0,b-a,n,p)$ 
we may restrict ourselves to the case $a=0$.
We will construct a Lie algebra isomorphism $\sigma$ from the extension of $S$ by $\langle D\rangle$
to the extension of $L$ by $\langle D\rangle$, where
$D=\bar y\frac{\partial}{\partial x}+(1+e)\frac{\partial}{\partial y}$.
(Using the same symbol for the derivation $D$ of both algebras should create no confusion.)
If we set
\begin{equation*}
\sigma(u_{ip^{n_2}+j})=-
x^{(p^{n_1}-i)}
y^{(p^{n_2}-j)}
\end{equation*}
for
$0<i\leq p^{n_1}$ and
$0<j\leq p^{n_2}$
with $(i,j)\not=(p^{n_1},p^{n_2})$,
we clearly have $D(\sigma(u_k))=\sigma(u_{k+1})=\sigma(Du_k)$ for all $k$.
In order to conclude that the bijective linear map determined by $\sigma$ is a Lie algebra isomorphism
it suffices to check that the defining relations of $S$ are satisfied in $H(2:\n;\omega_2)$.
In fact, we have
\begin{equation*}
 \begin{cases}{}
   [\sigma(u_{1+p^{n_2}}), \sigma(u_1)] = \{\bar x\bar y,x\}= \bar x
   y^{(p^{n_2-2})}= -\sigma(u_{2+p^{n_2}}), \\{}
   [\sigma(u_{1+p^{n_2}}), \sigma(u_{p^{n_2}})] = \{\bar x\bar y,y\}=
   -x^{(p^{n_1-2})}\bar y= \sigma(u_{1+2p^{n_2}}), \\{}
   [\sigma(u_{1+p^{n_2}}), \sigma(u_j)] =0 \qquad\text{in all other
     cases,}
 \end{cases}
\end{equation*}
where the last formula holds because $\{\bar x\bar y,x^{(k)} y^{(l)}\}=0$
unless $k,l\leq 1$, and also $\{\bar x\bar y,xy\}= \bar x
y^{(p^{n_2}-2)} y- x^{(p^{n_1}-2)} x \bar y=0$.

Incidentally, the sequence of the constituent lengths in the graded
Lie algebra of maximal class $AFS(0,b,n,p)$, which was computed in~\cite{CMN}
and utilized in~\cite{CN}, can now be more easily deduced from the isomorphism $\sigma$.
In fact, viewing the two-step centralizers in the corresponding (twisted) loop algebra of
$H(2:\n;\omega_2)$ and setting $q=p^{n_2}$
we see at once that
all components $L_k$ are centralized by $Y=(D+u_1)\otimes t=(D+x)\otimes t$, except when
$k\equiv 1 \pmod{q-1}$ and $q\not\equiv 1\pmod{p^n-1}$,
in which cases $L_k$ is centralized by $X=u_1\otimes t=x\otimes t$
(in the usual notation of~\cite{CMN} but with capital letters instead).

Another consequence of the isomorphism $\sigma$ is a formula
for carrying out explicit computations in the
Lie algebra of maximal class $AFS(0,b,n,p)$.
In fact, setting
\[
[u_{iq+j},u_{kq+l}]=
c(i,j,k,l)\cdot
u_{(i+k)q+(j+l)}
\]
for
$0<i,k\leq p^{n_1}$ and
$0<j,l\leq p^{n_2}$,
with $(i,j)$ and $(k,l)\not=(p^{n_1},p^{n_2})$,
we have
\begin{align*}
c(i,j,k,l)=&
-\binom{2p^{n_1}-i-k-1}{p^{n_1}-i}
\binom{2p^{n_2}-j-l-1}{p^{n_2}-j-1}
\\&
+
\binom{2p^{n_1}-i-k-1}{p^{n_1}-i-1}
\binom{2p^{n_2}-j-l-1}{p^{n_2}-j}.
\end{align*}
This rather unpleasant formula can be put into the slightly simpler form
\[
(-1)^{i+j}\cdot c(i,j,k,l)=
\binom{k}{p^{n_1}-i}
\binom{l-1}{p^{n_2}-j-1}
-
\binom{k-1}{p^{n_1}-i-1}
\binom{l}{p^{n_2}-j}
\]
by means of standard binomial coefficient manipulations
and the less standard but easily proved fact that the value of
$\binom{n}{k}$ modulo $p$ (for $n,k\in\Z$) is periodic in $n$ with period the smallest
power of $p$ which is greater than $k$.
Since any nonzero homogeneous element of weight $k\ge 2$ of $AFS(0,b,n,p)$
(realized as a loop algebra of $H(2:(n-n_2,n_2);\omega_2)$ as above)
can be uniquely written as a scalar multiple of
\[
[Y,X,Z_3,Z_4,\ldots,Z_k]=
u_k\otimes t
\]
with $Z_i\in\{X,-Y\}$,
the above formula allows one to multiply homogeneous elements
in $AFS(0,b,n,p)$ without any need for commutator expansions.
Finally, the method of iterated deflation described in~\cite{CMN} reduces computations in $AFS(a,b,n,p)$ to
computations in $AFS(0,b-a,n,p)$.

The isomorphism which we have seen above between the algebra of Albert and Frank $AF(0,b,n,p)$
and the Hamiltonian algebra
$H(2:\n;\omega_2)$ is perhaps better expressed in terms of the
basis $\{e_{\xi}\}$ of $AF(0,b,n,p)$.
As we have mentioned above, the general case of $AF(a,b,n,p)$ can be reduced to that of $AF(0,b-a,n,p)$
as in~\cite{CMN}, and we record the more general case in the following theorem
(which can, of course, also be proved by direct but rather tedious computation).

\begin{theorem}\label{thm:AF=Hamiltonian}
The algebra of Albert and Frank $AF(a,b,n,p)$ with Lie bracket given by~\eqref{eq:AF}
(in arbitrary prime characteristic $p$)
is isomorphic with the Hamiltonian algebra
$H(2:\n;\omega_2)$, where $\n=(n-b+a,b-a)$.
An isomorphism is given by the linear map
$\sigma:
AF(a,b,n,p)\to
H(2:\n;\omega_2)$
defined by the formula
\[
\sigma(e_{\xi})=
\sum
\xi^{ip^b+jp^a} x^{(i)}y^{(j)},
\]
where the summation is over all pairs
$(i,j)$ with $0\le i<n-b+a$, $0\le j<b-a$ and $(i,j)\not=(0,0)$.
The inverse map is given by
\[
\sigma^{-1}(x^{(i)}y^{(j)})=
-\sum_{\xi\in\F_{p^n}^\ast} \xi^{-ip^b-jp^a} e_\xi.
\]
\end{theorem}

This result gives an explicit realization of $H(2:\n;\omega_2)$
as a Block algebra.
Note that the subalgebra of $AF(a,b,n,p)$ which corresponds to the soluble subalgebra
$\langle x^{(i)}y^{(j)}:i+j>0\rangle$
of $H(2:\n;\omega_2)$ under the isomorphism $\sigma$
consists of all elements
$\sum_\xi c_\xi\,e_\xi$
with
$\sum_\xi c_\xi\,\xi^{p^a}=\sum_\xi c_\xi\,\xi^{p^b}=0$.
Because of this isomorphism (or by direct verification) this subalgebra of $AF(a,b,n,p)$
is maximal whenever the latter is simple
(that is, except when $p=2$ and either $b-a=1$ or $n-b+a=1$).
It is known that for $p>3$ this is the only maximal subalgebra of codimension two,
see~\cite[Theorem~2.16]{BO}.

Now we compute the derivation algebra of an algebra of Albert and Frank.
According to the identification of the algebras of Albert and Frank with suitable Hamiltonian algebras
given in Theorem~\ref{thm:AF=Hamiltonian}, we may consider their derivation algebras as known from the
general results on algebras of Cartan type quoted in Section~\ref{sec:cohomology}, but only for $p>3$.
Also, the derivation algebras of all Block algebras (thus including the algebras of Albert and Frank)
were computed already in~\cite{Bl}, again under the assumption that $p>3$.
The reason for this assumption in~\cite{Bl} was that genuine exceptions occur for $p=2,3$
in the more general case of Block algebras.
In particular, additional derivations may occur for algebras $H(2:\n;\omega_0)$ in low characteristics,
as we have illustrated in Remarks~\ref{rem:char3} and~\ref{rem:char2-derivations}.
Such exceptions do not arise for the algebras of Albert and Frank, that is, for Hamiltonian algebras
$H(2:\n;\omega_2)$.
However, Block's proof in~\cite{Bl} (in particular, Lemma~12) was not devised to deal with low characteristics
in this special case,
and we offer a variation of it which works in every characteristic.

Every specialization of the defining $\F_{p^n}$-grading of
the algebra of Albert and Frank $AF(a,b,n,p)$ to an $\F_p$-grading
(this terminology was explained in Section~\ref{sec:gradings})
gives rise to a derivation which multiplies each basis element
$e_{\xi}$ by its $\F_p$-degree.
These derivations turn out to be a basis for the space of outer derivations of $L=AF(a,b,n,p)$.
(In fact, these derivations span a maximal torus $T$ in $\Der(L)$,
and the given grading is the corresponding decomposition of $L$ into root spaces.)

\begin{theorem}\label{thm:AF-derivations}
The outer derivation algebra $\Der(L)/\ad(L)$ of the algebra of Albert and Frank
$AF(a,b,n,p)$ has dimension $n$.
More precisely, every derivation of degree zero 
with respect to the $\F_{p^n}$-grading of $L$ acts as
$D_{\pi}e_\xi=c_{\xi}\cdot e_\xi$, for some
additive map $\pi:\xi\mapsto c_\xi$ of $\F_{p^n}$ to itself;
all derivations of nonzero degree are inner.
Thus, a basis for the space of outer derivations is given by the derivations
$(D_{\mathrm{id}})^{p^s}e_\xi=\xi^{p^s}\cdot e_\xi$,
for $0\le s<n$.
\end{theorem}

\begin{proof}
Clearly $\Der(L)$ inherits an $\F_{p^n}$-grading from $L$, and the maps described are derivations of degree zero.
Conversely, let $D$ be a derivation of degree zero.
Then we have $De_{\xi}=c_{\xi}\cdot e_{\xi}$ for
$\xi\in\F_{p^n}^{\ast}$, for suitable scalars $c_{\xi}\in\F_{p^n}$, and we may set $c_{0}=0$.
The Leibniz rule $D[e_\xi,e_\eta]=[De_\xi,e_\eta]+[e_\xi,De_\eta]$
for $\xi,\eta\in\F_{p^n}^\ast$ yields that
$c_{\xi+\eta}=c_{\xi}+c_{\eta}$ provided $[e_{\xi},e_{\eta}]\not=0$.
The latter amounts to $\xi^{p^a}\eta^{p^b}-\xi^{p^b}\eta^{p^a}\not=0$,
which is satisfied as long as $\xi$ and $\eta$ do not span the same $\F_{p^c}$-subspace of $\F_{p^n}$,
where $c=(b-a,n)$.
But if $\xi$ and $\eta$ do span the same $\F_{p^c}$-subspace of $\F_{p^n}$,
we may choose $\theta$ outside this subspace and obtain that
$c_{\xi+\eta}+c_{\theta}=
c_{\xi+\eta+\theta}=
c_{\xi+\theta}+c_{\eta}=
c_{\xi}+c_{\eta}+c_{\theta}$
(which also settles the case where $\xi+\eta=0$).
Thus, the map $c:\F_{p^n}\to\F_{p^n}$ is additive.
These maps form an $n$-dimensional $\F_{p^n}$-space.

To simplify the calculations which follow we assume $a=0$, as we may, and set $q=p^b$.
Let $D$ be a derivation of degree $\tau\in\F_{p^n}^\ast$, hence
$De_{\xi}=c_{\xi}\cdot e_{\xi+\tau}$ for $\xi\in\F_{p^n}^\ast$
(recalling our convention that $e_0=0$),
where $c_\xi\in\F_{p^n}$, and we may set $c_{0}=c_{-\tau}=0$.
Then the condition
$D[e_\xi,e_\eta]=[De_\xi,e_\eta]+[e_\xi,De_\eta]$
translates into the equation
\[
(\xi\eta^q-\xi^q\eta)\cdot c_{\xi+\eta}=
((\xi+\tau)\eta^q-(\xi+\tau)^q\eta)\cdot c_{\xi}+
(\xi(\eta+\tau)^q-\xi^q(\eta+\tau))\cdot c_{\eta}
\]
for $\xi,\eta\in\F_{p^n}^{\ast}$ 
with $\xi+\eta\not=-\tau$.
However, this equation is trivially satisfied also if $\xi=0$ or $\eta=0$ or $\xi+\eta=-\tau$,
hence it holds for all $\xi,\eta\in\F_{p^n}$.

For the sake of clarity we only solve the equation in the case where $\tau=1$.
This is no loss because in general the equation can be reduced to this case by setting
$c'_{\xi}=c_{\xi\tau}$.
The condition for $D$ being a derivation can be written as
\begin{equation}\label{eq:Block}
(\xi\eta^q-\xi^q\eta)(c_{\xi+\eta}-c_{\xi}-c_{\eta})=
(\eta^q-\eta)\cdot c_{\xi}+
(\xi-\xi^q)\cdot c_{\eta}
\quad\textrm{for $\xi,\eta\in\F_{p^n}$.}
\end{equation}
One solution is $c_{\xi}=\xi^q-\xi$, which corresponds to the inner derivation $\ad e_1$.
We will prove that this is the only solution, up to scalars.

Setting $\eta=-\xi$ in equation~\eqref{eq:Block} we obtain that
$(\xi^q-\xi)(c_\xi+c_{-\xi})=0$, hence
$c_{-\xi}=-c_{\xi}$ holds provided $\xi^q-\xi\not=0$.
Our next goal is to show that the left-hand side of~\eqref{eq:Block} always vanishes.
We may assume that $\xi^q-\xi$ and $\eta^q-\eta$ do not both vanish,
otherwise the coefficient $\xi\eta^q-\xi^q\eta$ vanishes, too.
Then the left-hand side of~\eqref{eq:Block} can also be written as
$-(\xi\eta^q-\xi^q\eta)(c_{-\xi-\eta}+c_{\xi}+c_{\eta})$
if $\xi^q+\eta^q-\xi-\eta\not=0$, and as
$(\xi\eta^q-\xi^q\eta)(c_{\xi+\eta}+c_{-\xi}+c_{-\eta})$
otherwise.
Either expression is invariant under the cyclic substitution
$\xi\mapsto\eta\mapsto -\xi-\eta\mapsto\xi$.
Hence the right-hand side of~\eqref{eq:Block} is invariant, too, implying that
\begin{align*}
(\eta^q-\eta)\cdot c_{\xi}-
(\xi^q-\xi)\cdot c_{\eta}
&=
(-\xi^q-\eta^q+\xi+\eta)\cdot c_{\eta}-
(\eta^q-\eta)\cdot c_{-\xi-\eta}
\\
&=
(\xi^q-\xi)\cdot c_{-\xi-\eta}-
(-\xi^q-\eta^q+\xi+\eta)\cdot c_{\xi}
\end{align*}
or, equivalently,
\begin{align*}
(\xi^q-\xi)(c_{-\xi-\eta}+c_{\xi}+c_{\eta})=
(\eta^q-\eta)(c_{-\xi-\eta}+c_{\xi}+c_{\eta})=0.
\end{align*}
(A more conceptual view of the computation is the following:
invariance of the right-hand side of~\eqref{eq:Block}
under the substitution implies that
the (standard) vector product of the two vectors
$(c_\xi,c_\eta,c_{-\xi-\eta})$ and
$(\xi^q-\xi,\eta^q-\eta,-\xi^q-\eta^q+\xi+\eta)$
is a multiple of $(1,1,1)$;
since the vector product is orthogonal to both vectors we have
$c_{-\xi-\eta}+c_{\xi}+c_{\eta}=0$, provided the second vector is nonzero.)
Consequently, both sides of~\eqref{eq:Block} vanish for all $\xi,\eta\in\F_{p^n}$.
Hence
$(\eta^q-\eta)\cdot c_\xi=(\xi^q-\xi)\cdot c_\eta$,
and fixing $\eta$ such that $\eta^q-\eta\not=0$ we conclude that
$c_\xi=((\eta^q-\eta)^{-1} c_\eta)\cdot (\xi^q-\xi)$ for all $\xi$.
\end{proof}

We have already used Theorem~\ref{thm:AF-derivations} in Section~\ref{sec:cohomology}
to extend the scope of Theorem~\ref{thm:cocycles} to characteristic three.
We mention that a convenient basis for the second cohomology group of $AF(a,b,n,p)$
is obtained (in odd characteristic) from the outer derivations given in Theorem~\ref{thm:AF-derivations}
by the method explained in Section~\ref{sec:cohomology},
and consists of the cocycles
$\varphi((D_{\mathrm{id}})^{p^s})
(e_\xi,e_\eta)=\delta(\xi+\eta,0)\cdot\xi^{p^s}$
for $0\le s<n$.
In characteristic two these functions are not cocycles because they are not alternating and, in fact,
the second cohomology group vanishes, see Remark~\ref{rem:char2-cocycles}.

\section{Another thin grading of $H(2:\n;\omega_2)$}\label{sec:thin}

In this section we introduce a $\Z/(p^{n_2}-1)\Z \times (\Z/p\Z)^{n_1}$-grading of
$L=H(2:\n;\omega_2)$.
In the special case $n_1=1$, the grading will be (cyclic and) thin,
and the corresponding loop algebra will be the thin Lie algebra with all the diamonds
of finite type constructed in~\cite{CaMa:thin}.
In fact, the construction of the present grading is guided by
the structure of that thin Lie algebra.
Unlike in~\cite{CaMa:thin}, here we put no restriction on the (positive) characteristic.

To briefly outline the construction of the grading, we first note that such a grading gives rise
to a $\Z/(p^{n_2}-1)$-grading and to a $(\Z/p\Z)^{n_1}$-grading by specialization.
Conversely, the original grading can be recovered from these two specializations.
We start with the $\Z/(p^{n_2}-1)\Z$-grading of $L$ mentioned at the beginning of
Section~\ref{sec:gradings} where $(R,S)=(0,-1)$.
Its components are
\begin{equation*}
\bar L_{1}=
\langle
x^{(i)}
\mid
i=1,\ldots,p^{n_1}-1
\rangle
+
\langle
x^{(i)}\bar y
\mid
i=0,\ldots,p^{n_1}-1
\rangle,
\end{equation*}
and
\begin{equation*}
\bar L_{1-j}=
\langle
x^{(i)}y^{(j)}
\mid
i=0,\ldots,p^{n_1}-1
\rangle
\end{equation*}
for $j=1,\ldots,p^{n_2}-2$.
The component of degree
zero is isomorphic with a Zassenhaus algebra $W(1:n_1)$.
As we saw in Section~\ref{sec:Cartan}, the latter has a grading over the additive group of the field
$\F_{p^{n_1}}$, which is simply the root space decomposition of $W(1:n_1)$
with respect to an appropriate one-dimensional Cartan subalgebra, say $\langle e_0\rangle$.
Viewing $W(1:n_1)$ as a subalgebra of $L$, the decomposition of $L$ into weight spaces
with respect to $\ad\,e_0$ extends the grading of $W(1:n_1)$ to a grading of $L$,
which turns out to be over the same group $\F_{p^{n_1}}$.
Since $e_0\in\bar L_0$, it normalizes every component $\bar L_j$, and it follows that
this $\F_{p^{n_1}}$-grading of $L$ together with the $\Z/(p^{n_2}-1)\Z$-grading yield the
desired $\Z/(p^{n_2}-1)\Z \times \F_{p^{n_1}}$-grading.

As announced above, we identify the subalgebra
\begin{equation*}
W=\langle
x^{(i)}y\mid
i=0,\ldots,p^{n_1}-1
\rangle
\end{equation*}
of $L$ with the Zassenhaus algebra $W(1:n_1)$ via $E_i=x^{(i+1)}y$, for
$i=-1,\ldots,p^{n_1}-2$.
This is justified by
\begin{equation*}
\{x^{(i)}y,
x^{(k)}y\}
=
\left(
{i+k-1 \choose k-1}-
{i+k-1 \choose i-1}
\right)
x^{(i+k-1)}y.
\end{equation*}
Recall from Section~\ref{sec:Cartan} that an $\F_{p^{n_1}}$-graded basis
of $W$ is given by
\begin{equation*}
\begin{cases}
  e_0= y+\bar x y
  &\\
  e_{\alpha}= \bar x y+ \sum_{i=0}^{p^{n_1}-1} \alpha^{i} x^{(i)}y &
  \text{ for } \alpha\in\F_{p^n}^{\ast}
\end{cases}
\end{equation*}
(where the former formula can be seen as a special case of
the latter, but we have kept them separated for the sake of clarity, here and below).

It is a simple matter to find the weight spaces of $\ad e_0$ on each component $\bar L_j$,
starting from the formula
\begin{align*}
  \{e_0, x^{(i)} y^{(j)}\} &= \{y+\bar x y, x^{(i)}
  y^{(j)}\}=\\
  & =
\begin{cases}
  -j x^{(p^{n_1}-2)} y^{(j)}
  &\text{if $i=0$,}\\
  (1+e) y^{(j)}+ (1+j) \bar x y^{(j)}
  &\text{if $i=1$,}\\
  x^{(i-1)} y^{(j)} &\text{if $i>1$.}
\end{cases}
\end{align*}
In particular,
eigenvectors for $\ad\, e_0$ on $\bar L_{1-j}$, for
$j=1,\ldots,p^{n_2}-1$, are
\begin{equation*}
\begin{cases}
  e_{1-j,0}= y^{(j)}+ j\bar x y^{(j)}
  &\\
\displaystyle
  e_{1-j,\alpha}= j\bar x y^{(j)}+ \sum_{i=0}^{p^{n_1}-1} \alpha^{i}
  x^{(i)}y^{(j)} & \text{for $\alpha\in\F_{p^{n_1}}^{\ast}$},
\end{cases}
\end{equation*}
extending the complete set of eigenvectors for $\ad\, e_0$ on $W$
given by its basis elements $e_{0,\alpha}=e_{\alpha}$.
Our notation is chosen so that
$e_{1-j,\alpha}\in\bar L_{1-j}$, and
\begin{equation*}
\{e_0,e_{1-j,\alpha}\}=
\alpha e_{1-j,\alpha}.
\end{equation*}
Also, the subscript $1-j$ in $e_{1-j,\alpha}$ will be read
modulo $p^{n_2}-1$, so that, for example, $e_{1,\alpha}$ will be the same as
$e_{2-p^{n_2},\alpha}$.
However, beware that $e_{1,\alpha}$ is not
what one would obtain by putting $j=0$ in the
formulas.

The formulas above give complete sets of eigenvectors for $\ad\, e_0$
on $\bar L_{1-j}$ only for $j=1,\ldots,p^{n_2}-2$.
In fact, $\ad\, e_0$ acts on $\bar L_{2-p^{n_2}}$ as a sum of two-dimensional Jordan blocks,
say $\langle e_{1,\alpha}, \bar e_{1,\alpha} \rangle$, one for each eigenvalue $\alpha\in\F_{p^{n_1}}^{\ast}$.
(These will be the diamonds in the thin grading, in case $n_1=1$.)
It is convenient to assume that
\begin{equation*}
\{e_0,\bar e_{1,\alpha}\}=
\alpha \bar e_{1,\alpha}-
\alpha e_{1,\alpha},
\end{equation*}
which we may because $\alpha\not=0$.
This determines $\bar e_{1,\alpha}\in \bar L_1$ modulo
$\langle e_{1,\alpha}\rangle$, and we choose to set
\begin{equation*}
\bar e_{1,\alpha}=
\sum_{i=1}^{p^{n_1}-1}
\alpha^{i}
(x^{(i)}-
i x^{(i)}\bar y)
\end{equation*}
for $\alpha\in\F_{p^{n_1}}^{\ast}$.
It will be convenient to allow also the value $\alpha=0$
in the above formula, and thus set $\bar e_{1,0}=0$.
(In case $n_1=1$ the component
$\langle e_{1,0}, \bar e_{1,0} \rangle=\langle e_{1,0}\rangle$
corresponds to a {\em fake} diamond in the thin grading.)

In the following lemma we summarize what we have obtained so far, and we include
the inverse formulas which give the
divided powers in terms of the elements $e_{j,\alpha}$ and $\bar e_{j,\alpha}$.

\begin{lemma}\label{lemma:thin}
The Hamiltonian algebra $L=H(2:\n;\omega_2)$ admits a
$\Z/(p^{n_2}-1)\Z \times \F_{p^{n_1}}$-grading given by $L=\bigoplus L_{k,\alpha}$, where
$L_{k,\alpha}= \langle e_{k,\alpha} \rangle$, or
$L_{k,\alpha}= \langle e_{k,\alpha}, \bar e_{k,\alpha} \rangle$ whenever the latter makes sense, and
the basis elements involved are given by
\begin{equation*}
\begin{cases}
\displaystyle
  e_{1-j,\alpha}=
  j\bar x y^{(j)}+
  \sum_{i=0}^{p^{n_1}-1} \alpha^{i}
  x^{(i)}y^{(j)}
  &\text{for $\alpha\in\F_{p^{n_1}}$ and $j=1,\ldots,p^{n_2}-1$},\\
\displaystyle
  \bar e_{1,\alpha}=
  \sum_{i=1}^{p^{n_1}-1}
  \alpha^{i}
  (x^{(i)}-
  i x^{(i)}\bar y)
  &\text{for $\alpha\in\F_{p^{n_1}}^{\ast}$},
\end{cases}
\end{equation*}
where $0^0$ should be understood as $1$.
Conversely, for
$i=1,\ldots,p^{n_1}-1$ and $j=1,\ldots,p^{n_2}-1$ we have
\begin{equation*}
\begin{cases}\displaystyle
  x^{(i)}y^{(j)}= -\sum_{\alpha\in\F_{p^{n_1}}} \alpha^{p^{n_1}-1-i}
  e_{1-j,\alpha},\\\displaystyle
  y^{(j)}= e_{1-j,0}+ j\sum_{\alpha\in\F_{p^{n_1}}}
  e_{1-j,\alpha},\\\displaystyle
  x^{(i)}= -\sum_{\alpha\in\F_{p^{n_1}}} \alpha^{p^{n_1}-1-i} (i
  e_{1,\alpha}+\bar e_{1,\alpha}),
\end{cases}
\end{equation*}
where in the last formula it is understood that $\bar e_{1,0}=0$.
\end{lemma}

For later use we record the products of the basis elements given in Lemma~\ref{lemma:thin}. We have
\begin{align}
&\label{eq:mult1}
\{e_{1-j,\alpha},e_{1-l,\beta}\}=
\left(\beta\binom{j+l-1}{l}-\alpha\binom{j+l-1}{j}\right)
e_{2-j-l,\alpha+\beta}
\intertext{for $j,l=1,\ldots,p^{n_2}-1$, provided
$1\le j+l-1\le p^{n_2}-1$,
and zero otherwise,}
&\label{eq:mult2}
\{e_{1-j,\alpha},\bar e_{1,\beta}\}=
\beta
e_{2-j,\alpha+\beta}
\quad\text{for $j=2,\ldots,p^{n_2}-1$},
\\&\label{eq:mult3}
\{e_{0,\alpha},\bar e_{1,\beta}\}=
-\beta
e_{1,\alpha+\beta}+\beta\bar e_{1,\alpha+\beta},
\\&\label{eq:mult4}
\{\bar e_{1,\alpha},\bar e_{1,\beta}\}=0.
\end{align}
It is not necessary to carry out the computations in full in order to prove these formulas,
if we use the fact that the basis is graded, according to Lemma~\ref{lemma:thin}.
For example, the result of the product in formula~\eqref{eq:mult1}
must be a scalar multiple of $e_{2-j-l,\alpha+\beta}$, and we only have to find that scalar.
Since the monomial $y^{(j)}$ appears in $e_{1-j,\alpha}$ with coefficient $1$,
it is enough to compute the coefficient of the monomial $y^{(j+l-1)}$ in the product $\{e_{1-j,\alpha},e_{1-l,\beta}\}$.
The products of the only relevant terms, namely
\[
\{y^{(j)},\beta xy^{(l)}\}+\{\alpha xy^{(j)},y^{(l)}\}=
\left(\beta\binom{j+l-1}{l}-\alpha\binom{j+l-1}{j}\right)
y^{(j+l-1)},
\]
yield the desired conclusion.
Formula~\eqref{eq:mult2} follows similarly by computing
$
\{y^{(j)},\beta x\}=\beta y^{(j-1)}.
$
Formula~\eqref{eq:mult4} follows because the monomial $y^{(p^{n_2}-2)}$ does not appear in the result.
The proof of formula~\eqref{eq:mult3} is just slightly more involved because
$\{e_{0,\alpha},\bar e_{1,\beta}\}$ must be a linear combination of
$e_{1,\alpha+\beta}$  and $\bar e_{1,\alpha+\beta}$.
Since the coefficients of $\bar y$ and $x$ are $1$ and $0$ in the former,
and $0$ and $\alpha+\beta$ in the latter,
the conclusion follows by computing the only relevant products, which are
$\{y,-\beta x\bar y\}=-\beta\bar y$ and
\[
\{y,\beta^2 x^{(2)}\}+
\{\alpha xy,\beta x\}=
\beta(\alpha+\beta)x.
\]
(For the dubious case where $\alpha+\beta=0$ recall from Lemma~\ref{lemma:thin} that we have set $\bar e_{1,0}=0$.)

Now we restrict our attention to the case of main interest for us
by letting $n_1=1$, and introducing the shorthand $q=p^{n_2}$, so that
$L$ has dimension $pq-1$.
Then the subspaces $L_{k,\alpha}$ form a
grading of $L$ over the cyclic group $\Z/(q-1)\Z \times \F_p$.
Choosing $(1,1)$ as a generator of the latter,
we will show that the grading is thin, according to Definition~\ref{def:thin-fd},
and that all diamonds are of finite type.

The diamonds occur in all degrees congruent to $1$ modulo $q-1$, with the exception of degree $q$.
We set
\[
X=\bar e_{1,1}=
\sum_{i=1}^{p-1}
(x^{(i)}-
i x^{(i)}\bar y)
\quad\text{and}\quad
Y=e_{1,1}=
-\bar x\bar y+
\sum_{i=0}^{p-1}
x^{(i)}\bar y,
\]
and check that the covering property is satisfied.
Since
\begin{align*}
  &
  \{\bar e_{1,\alpha},X\}=0,\\
  & \{e_{0,\alpha},X\}= \bar e_{1,\alpha+1} -e_{1,\alpha+1}
  \quad\text{(where $\bar e_{1,0}=0$)},\\
  & \{e_{1-j,\alpha},X\}= e_{2-j,\alpha+1}
  \quad\text{for $j=2,\ldots q-1$},\\
  & \{\bar e_{1,\alpha},Y\}=
  -\alpha e_{2,\alpha+1},\\
  & \{e_{0,\alpha},Y\}=
  (\alpha+1) e_{1,\alpha+1},\\
  & \{e_{1-j,\alpha},Y\}=0 \quad\text{for $j=2,\ldots q-1$},
\end{align*}
we have
\begin{equation*}
\{L_{k,\alpha},\langle X,Y\rangle\}=
L_{k+1,\alpha+1}
\quad\textrm{for $k\not\equiv 1 \pmod{(q-1)}$,}
\end{equation*}
which shows that the covering property holds in these degrees, since $L_{k,\alpha}$ is one-dimensional.
To see that the loop algebra is thin, it remains to check the covering
property on the two-dimensional components.  At the same time we will
check the diamond types (including the {\em fake} diamond $L_{1,0}$, which is in degree $q$).
Since
\begin{align*}
  & \{e_{0,\alpha},X,X\}= \{-e_{1,\alpha+1},X\}=
  -e_{2,\alpha+2},\\
  & \{e_{0,\alpha},X,Y\}= \{\bar e_{1,\alpha+1},Y\}=
  -(\alpha+1)e_{2,\alpha+2},\\
  & \{e_{0,\alpha},Y,X\}= \{(\alpha+1)e_{1,\alpha+1},X\}=
  (\alpha+1)e_{2,\alpha+2},\\
  & \{e_{0,\alpha},Y,Y\}= 0,
\end{align*}
if $a,b$ are scalars we have
$\{e_{0,\alpha},aX+bY,\langle X,Y\rangle\}=L_{2,\alpha+2}$
unless $a=b=0$ in case $\alpha\not=-1$, and unless $a=0$ in case $\alpha=-1$
(the case of the {\em fake} diamond).
It follows that the grading is thin.
Furthermore, the elements just above the (possibly {\em fake}) diamonds are
those of the form $V_{\alpha}= e_{0,-\alpha-1}$ and satisfy
\begin{equation*}
\{V_{\alpha},Y,X\}=
\alpha \{V_{\alpha},X,X\},
\end{equation*}
hence the diamond
$\langle\{V_{\alpha},X\},\{V_{\alpha},Y\}\rangle$ in degree $q+\alpha(q-1)$
(which is {\em fake} for $\alpha=0$) has type $\alpha$.
Note that most of the computations done in this paragraph need not be carried out in full,
provided we use the fact that the elements
$e_{k,\alpha}$ and $\bar e_{1,\alpha}$ form a graded basis of $L$, according to Lemma~\ref{lemma:thin}.

We have completed the proof of the following result.

\begin{theorem}
In case $n_1=1$ the cyclic grading of $L=H(2:\n;\omega_2)$ defined in Lemma~\ref{lemma:thin}
is a thin grading (with respect to the generator $(1,1)$ of the grading group).
The diamond types are all finite, and follow an arithmetic progression.
\end{theorem}

Note that
$\{L_{k,\alpha},Y\}=0$
for $k\not\equiv 0,1 \pmod{(q-1)}$,
and also
$\{L_{0,-1},Y\}=\{L_{1,0},Y\}=0$;
hence $Y$ centralizes all one-dimensional components which do not immediately precede a diamond.

We remark that thin Lie algebras with diamonds of both finite and infinite types are constructed in~\cite{AviMat:-1}.
The finite diamond types there follow an arithmetic progression, but they are separated by sequences of
constant length of diamonds of infinite type.

\section{Another realization of $H(2:(1,n_2);\omega_2)$ as a Block algebra}\label{sec:Andrea}

In~\cite{CaMa:thin} we constructed a simple Lie algebra $A$ with a thin grading
and all diamonds of finite type.
Since the location and types of the diamonds of $A$ match those of
the algebra $L=H(2:(1,n_2);\omega_2)$ with the grading defined in the preceding section
it follows that the loop algebras of $L$ and $A$ are isomorphic.
Actually, this yields an isomorphism between $L$ and $A$ themselves, but it may be necessary
to extend the ground field for this.
The purpose of this section is to justify these claims.

We recall the definition of $A$ from~\cite{CaMa:thin}, with slight notational changes
to avoid clash with the notation of the present paper.
Let $p$ be any prime, let $q=p^{n_2}$, and set $S^{\ast}=(\F_p\times\F_q)\setminus\{(0,0)\}$.
The $\F_{q}$-vector space $A$ with basis
$\{f_{u,\alpha}\mid (u,\alpha)\in S^{\ast
}\}$
becomes a Lie algebra over $\F_{q}$ by defining
\[
 [f_{u,\alpha}, f_{v,\beta}] = (v\alpha-u\beta)\cdot f_{u+v,\alpha+\beta},
\]
where we read $0\cdot f_{0,0}$ as zero.
By construction $A$ is a Block algebra,
and so it is simple and isomorphic
to a Hamiltonian algebra $H(2:\n;\omega_2)$
(at least when $p>3$, and possibly after extending the ground field, see Section~\ref{sec:Hamiltonian}).
We will construct an isomorphism explicitly.

It was shown in~\cite{CaMa:thin} that $A$ has a thin grading
where the component of degree one is spanned by
\begin{equation*}
 x = f_{1, 0} \quad\textrm{and}\quad
 y = \sum_{\alpha \in \F_{q}} f_{1, \alpha}.
\end{equation*}
(The assumption $p>3$ stated in~\cite{CaMa:thin} was only needed in later sections, when considering presentations
for the thin algebra.)
An easy induction showed that
\begin{equation*}
[y,
\U{x}{j-1}]=
 \sum_{\alpha \in \F_{q}}
 \alpha^{j-1} f_{j, \alpha}
\qquad\text{for $j\ge 1$}
\end{equation*}
(where our convention that $0^0=1$ intervenes when $j=1$).
This element spans the component of degree $j$ of the grading,
except when $j\equiv 1\pmod{q-1}$ but $j\not\equiv 0\pmod{p}$,
in which cases the component is two-dimensional.
An element just above a diamond of type $\lambda-1$ is
\begin{equation*}
v_{\lambda-1}=[y,
\U{x}{\lambda(q-1)-1}]=
  \sum_{\alpha \in \F_{q}^{\ast}}
 \alpha^{-1} f_{-\lambda, \alpha},
\end{equation*}
whence
\begin{equation*}
[y,
\U{x}{\lambda(q-1)-1},y]=
 (\lambda - 1)
 \sum_{\gamma \in \F_{q}}
 f_{-\lambda+1, \gamma}.
\end{equation*}
In particular, we have
\begin{equation*}
\begin{cases}
\displaystyle
  {}[y, \U{x}{p(q-1)}]= \sum_{\alpha \in \F_{q}} \alpha^{(q - 1) p}
  f_{1 + (q - 1) p, \alpha} = \sum_{\alpha \in \F_{q}^{*}} f_{1,
    \alpha} = y - x,
  \\
\displaystyle
  {}[y, \U{x}{p(q-1)-1},y]= -\sum_{\gamma \in \F_{q}} f_{1, \gamma}
  =-y.
\end{cases}
\end{equation*}

Now we carry out analogous computations in $L=H(2:(1,n_2);\omega_2)$
with the thin grading constructed in the previous section, using the formulas
obtained there.
In particular, the component of degree one in the grading is spanned by
$X=\bar e_{1,1}$ and $Y=e_{1,1}$, and we find that
\begin{equation*}
\{Y,
\U{X}{j-1}\}=
(-1)^{\lfloor\frac{j-1}{q-1}\rfloor}
(e_{j,j}-
\delta_{q-1}(j,1)
\bar e_{j,j})
\qquad\text{for $j> 1$}
\end{equation*}
(but not for $j=1$),
where $\delta_{q-1}(s,t)$ equals one if $s\equiv t\pmod{q-1}$, and zero otherwise.
An element just above a diamond
of type $\lambda-1$ is
\begin{equation*}
V_{\lambda-1}=\{Y,
\U{X}{\lambda(q-1)-1}\}=
(-1)^{\lambda-1}
e_{0,-\lambda},
\end{equation*}
whence
\begin{equation*}
\{Y,
\U{X}{\lambda(q-1)-1},Y\}=
(-1)^{\lambda}(\lambda-1)
e_{1,-\lambda+1}.
\end{equation*}
In particular, we have
\begin{equation*}
\begin{cases}
  {}\{Y, \U{X}{p(q-1)}\}= -e_{1,1}+\bar e_{1,1}= -Y+X, \\
  {}\{Y, \U{X}{p(q-1)-1},Y\}= e_{1,1}=Y.
\end{cases}
\end{equation*}

The sign discrepancy of the latter formulas with the analogous ones for $A$ seen above
means nothing if we work with the loop algebras of $A$ and $L$ (which are isomorphic),
but suggests that $A$ may not be isomorphic with $L=H(2:(1,n_2);\omega_2)$ over $\F_q$.
We extend the ground field to $\F_{q^2}$ and fix $\varepsilon\in\F_{q^2}$ with
$\varepsilon^{q-1}=-1$.
A well-defined linear map $\tau: A \to L$ is obtained
by setting
$\tau(x)=\varepsilon X$,
$\tau(y)=\varepsilon Y$,
and
\begin{equation*}
\begin{cases}
  \tau([y,\U{x}{j-2},x]) = \varepsilon^j \{Y,\U{X}{j-2},X\},
  \\
  \tau([y,\U{x}{j-2},y]) = \varepsilon^j \{Y,\U{X}{j-2},Y\}
\end{cases}
\end{equation*}
for all $j\ge 2$.
The computations done in Section~\ref{sec:thin} and in~\cite{CaMa:thin} show that
$y$ and $Y$ centralize all one-dimensional components which do not immediately precede a diamond,
in the respective algebras, and also that the diamond types match.
It follows that
$\tau([w,z])=\{\tau(w),\tau(z)\}$
for all $w\in A$ and $z=x$ or $y$.
Since $x$ and $y$ generate $A$
we conclude that the map $\tau$ is an isomorphism of Lie algebras
(because the set of $z\in A$ such that
$\tau([w,z])=\{\tau(w),\tau(z)\}$ for all $w\in A$ is a subalgebra).

Finding explicit formulas for the isomorphism $\tau$ in terms of convenient bases
for $A$ and $L$ is a matter of some manipulations, and we only record the final result.

\begin{theorem}\label{thm:Andrea=Hamiltonian}
The Block algebra $A$ defined in~\cite{CaMa:thin} and described at the beginning of this section
is isomorphic with the Hamiltonian algebra
$H(2:(1,n_2);\omega_2)$.
An isomorphism is given by the linear map
$\tau:
A\to
H(2:(1,n_2);\omega_2)$
defined by the formulas
(both for $\alpha\in\F_p$, and where $\varepsilon^{q-1}=-1$)
\begin{equation*}
\begin{cases}
\displaystyle
  \tau(f_{\alpha, \beta}) = -\sum_{k=1}^{q-1} \beta^{-k+1}
  \varepsilon^k e_{k,\alpha} +\varepsilon \bar e_{1,\alpha}
  &\text{for $\beta\in\F_q^{\ast}$},
  \\
\displaystyle
  \tau(f_{\alpha,0}) = \varepsilon \bar e_{1,\alpha},
\end{cases}
\end{equation*}
in terms of the basis $\{e_{k,\alpha},\bar e_{1,\alpha}\}$ of $H(2:(1,n_2);\omega_2)$
defined in Lemma~\ref{lemma:thin},
or
\begin{equation*}
\begin{cases}
\displaystyle
  \tau(f_{\alpha, \beta}-f_{\alpha,0}) =
  \sum_{j=1}^{q-1} j\varepsilon^{1-j} \beta^{j} \bar xy^{(j)} +
  \sum_{i=0}^{p-1} \sum_{j=1}^{q-1} \varepsilon^{1-j}
  \alpha^i\beta^{j} x^{(i)}y^{(j)}
  &\text{for $\beta\in\F_q^{\ast}$},
  \\
\displaystyle
  \tau(f_{\alpha,0})=
  \varepsilon
  \sum_{i=1}^{p-1}
  \alpha^{i}(x^{(i)}-i x^{(i)}\bar y)
\end{cases}
\end{equation*}
in terms of divided powers.
\end{theorem}

Note that 
there appears to be a certain amount of symmetry between the formulas for $\tau(f_{0, \beta})$ and
$\tau(f_{\alpha,0})$, because
$
\tau(f_{0,\beta})
=
\sum_{j=1}^{q-1}
\varepsilon^{1-j}
\beta^{j}
(y^{(j)}+
j \bar xy^{(j)}).
$

The following cocycles of $A$ were introduced in~\cite{CaMa:thin}, and it was proved that they form
a basis for the second cohomology group of $A$, if $p>3$:
\[
 \varphi_{r} (f_{u, \alpha}, f_{v, \beta})
 =
 \begin{cases}
  \alpha^{p^{r}} & \text{if $\alpha + \beta = 0$
                         and $u + v = 0$}\\
  0 & \text{otherwise,}
 \end{cases}
\]
for $1 \le r \le n_2$, and
\[
 \psi (f_{u, \alpha}, f_{v, \beta})
 =
 \begin{cases}
  u & \text{if $u + v = 0$}\\
  0 & \text{otherwise.}
 \end{cases}
\]

In view of the isomorphism $\tau$ given in Theorem~\ref{thm:Andrea=Hamiltonian},
we now make the connection with the cocycles of $L=H(2:(1,n_2);\omega_2)$
described in Theorem~\ref{thm:cocycles}.
Call $\Phi_r$ and $\Psi$ the cocycles obtained by pulling the cocycles $\varphi_r$ and $\psi$
of $A$ back to $L$ via the isomorphism $\tau^{-1}$.
A straightforward computation for $\Psi$, and a slightly more involved one for $\Phi_r$,
both of which we omit, show that
$\Psi=-\varepsilon^{-2}\psi_1$
and that
\begin{equation*}
\Phi_r=
\varepsilon^{p^r-2}
\left(\varphi_r
+\varphi(
\ad xy^{(q-p^r)}
)\right),
\end{equation*}
in the notation of Section~\ref{sec:cohomology}.

\bibliography{References}

\newcommand{\etalchar}[1]{$^{#1}$}
\def\cprime{$'$} \def\cprime{$'$}
\providecommand{\bysame}{\leavevmode\hbox to3em{\hrulefill}\thinspace}
\providecommand{\MR}{\relax\ifhmode\unskip\space\fi MR }
\providecommand{\MRhref}[2]{%
  \href{http://www.ams.org/mathscinet-getitem?mr=#1}{#2}
}
\providecommand{\href}[2]{#2}
\begin{thebibliography}{BGO{\etalchar{+}}89}

\bibitem[AF55]{AF}
A.~A. Albert and M.~S. Frank, \emph{Simple {L}ie algebras of characteristic
  {$p$}}, Univ. e Politec. Torino. Rend. Sem. Mat. \textbf{14} (1954--55),
  117--139. \MR{MR0079222 (18,52a)}

\bibitem[AJ01]{AviJur}
M.~Avitabile and G.~Jurman, \emph{Diamonds in thin {L}ie algebras}, Boll.
  Unione Mat. Ital. Sez. B Artic. Ric. Mat. (8) \textbf{4} (2001), no.~3,
  597--608. \MR{MR1859998 (2003a:17038)}

\bibitem[AM]{AviMat:Nottingham}
M.~Avitabile and S.~Mattarei, \emph{Thin loop algebras of {A}lbert-{Z}assenhaus
  algebras}, submitted.

\bibitem[AM05]{AviMat:-1}
\bysame, \emph{Thin {L}ie algebras with diamonds of finite and infinite type},
  J. Algebra \textbf{293} (2005), no.~1, 34--64.

\bibitem[Avi99]{Avi:thesis}
M.~Avitabile, \emph{Some loop algebras of {Hamiltonian} {Lie} algebras}, Ph.D.
  thesis, Trento, November 1999.

\bibitem[Avi02]{Avi}
Marina Avitabile, \emph{Some loop algebras of {H}amiltonian {L}ie algebras},
  Internat. J. Algebra Comput. \textbf{12} (2002), no.~4, 535--567.
  \MR{MR1919687 (2003e:17013)}

\bibitem[BCS92]{BCS}
R.~Brandl, A.~Caranti, and C.~M. Scoppola, \emph{Metabelian thin {$p$}-groups},
  Quart. J. Math. Oxford Ser. (2) \textbf{43} (1992), no.~170, 157--173.
  \MR{MR1164620 (93d:20041)}

\bibitem[BGO{\etalchar{+}}89]{BGOSW}
G.~M. Benkart, T.~B. Gregory, J.~M. Osborn, H.~Strade, and R.~L. Wilson,
  \emph{Isomorphism classes of {H}amiltonian {L}ie algebras}, Lie algebras,
  Madison 1987, Lecture Notes in Math., vol. 1373, Springer, Berlin, 1989,
  pp.~42--57. \MR{MR1007323 (91e:17014)}

\bibitem[BKK95]{BKK}
Georgia Benkart, Alexei~I. Kostrikin, and Michael~I. Kuznetsov,
  \emph{Finite-dimensional simple {L}ie algebras with a nonsingular
  derivation}, J. Algebra \textbf{171} (1995), no.~3, 894--916. \MR{MR1315926
  (96b:17020)}

\bibitem[Blo58]{Bl}
Richard Block, \emph{New simple {L}ie algebras of prime characteristic}, Trans.
  Amer. Math. Soc. \textbf{89} (1958), 421--449. \MR{MR0100010 (20 \#6446)}

\bibitem[Blo69]{Bl:differentiably-simple}
Richard~E. Block, \emph{Determination of the differentiably simple rings with a
  minimal ideal.}, Ann. of Math. (2) \textbf{90} (1969), 433--459.
  \MR{MR0251088 (40 \#4319)}

\bibitem[BM86]{BenMoo}
G.~M. Benkart and R.~V. Moody, \emph{Derivations, central extensions, and
  affine {L}ie algebras}, Algebras Groups Geom. \textbf{3} (1986), no.~4,
  456--492. \MR{MR901810 (89a:17027)}

\bibitem[BO88]{BO}
Georgia Benkart and J.~Marshall Osborn, \emph{Toral rank one {L}ie algebras},
  J. Algebra \textbf{115} (1988), no.~1, 238--250. \MR{MR937612 (89d:17016)}

\bibitem[Bra88]{Br}
Rolf Brandl, \emph{The {D}ilworth number of subgroup lattices}, Arch. Math.
  (Basel) \textbf{50} (1988), no.~6, 502--510. \MR{MR948264 (89e:20054)}

\bibitem[BW82]{BlWil:rank-two}
Richard~E. Block and Robert~Lee Wilson, \emph{The simple {L}ie {$p$}-algebras
  of rank two}, Ann. of Math. (2) \textbf{115} (1982), no.~1, 93--168.
  \MR{MR644017 (83j:17008)}

\bibitem[Car97]{Car:Nottingham}
A.~Caranti, \emph{Presenting the graded {L}ie algebra associated to the
  {N}ottingham group}, J. Algebra \textbf{198} (1997), no.~1, 266--289.
  \MR{MR1482983 (99b:17019)}

\bibitem[Car98]{Car:thin_addendum}
\bysame, \emph{Thin groups of prime-power order and thin {L}ie algebras: an
  addendum}, Quart. J. Math. Oxford Ser. (2) \textbf{49} (1998), no.~196,
  445--450. \MR{MR1660053 (99k:17036)}

\bibitem[Car99]{Car:Zassenhaus-three}
\bysame, \emph{Loop algebras of {Z}assenhaus algebras in characteristic three},
  Israel J. Math. \textbf{110} (1999), 61--73. \MR{MR1750443 (2001d:17019)}

\bibitem[Car01]{Carrara}
Claretta Carrara, \emph{({F}inite) presentations of the
  {A}lbert-{F}rank-{S}halev {L}ie algebras}, Boll. Unione Mat. Ital. Sez. B
  Artic. Ric. Mat. (8) \textbf{4} (2001), no.~2, 391--427. \MR{MR1831996
  (2003d:17020)}

\bibitem[CJ99]{CaJu:quotients}
A.~Caranti and G.~Jurman, \emph{Quotients of maximal class of thin {L}ie
  algebras. {T}he odd characteristic case}, Comm. Algebra \textbf{27} (1999),
  no.~12, 5741--5748. \MR{MR1726275 (2001a:17042a)}

\bibitem[CM99]{CaMa:thin}
A.~Caranti and S.~Mattarei, \emph{Some thin {L}ie algebras related to
  {A}lbert-{F}rank algebras and algebras of maximal class}, J. Austral. Math.
  Soc. Ser. A \textbf{67} (1999), no.~2, 157--184, Group theory. \MR{MR1717411
  (2000j:17036)}

\bibitem[CM04]{CaMa:Nottingham}
\bysame, \emph{Nottingham {L}ie algebras with diamonds of finite type},
  Internat. J. Algebra Comput. \textbf{14} (2004), no.~1, 35--67. \MR{MR2041537
  (2004j:17027)}

\bibitem[CMN97]{CMN}
A.~Caranti, S.~Mattarei, and M.~F. Newman, \emph{Graded {L}ie algebras of
  maximal class}, Trans. Amer. Math. Soc. \textbf{349} (1997), no.~10,
  4021--4051. \MR{MR1443190 (98a:17027)}

\bibitem[CMNS96]{CMNS}
A.~Caranti, S.~Mattarei, M.~F. Newman, and C.~M. Scoppola, \emph{Thin groups of
  prime-power order and thin {L}ie algebras}, Quart. J. Math. Oxford Ser. (2)
  \textbf{47} (1996), no.~187, 279--296. \MR{MR1412556 (97h:20036)}

\bibitem[CN00]{CN}
A.~Caranti and M.~F. Newman, \emph{Graded {L}ie algebras of maximal class.
  {II}}, J. Algebra \textbf{229} (2000), no.~2, 750--784. \MR{MR1769297
  (2001g:17041)}

\bibitem[Dzh84]{Dzhu:central}
A.~S. Dzhumadil{\cprime}daev, \emph{Central extensions and invariant forms of
  {L}ie algebras of positive characteristic of {C}artan types}, Funktsional.
  Anal. i Prilozhen. \textbf{18} (1984), no.~4, 77--78. \MR{MR775937
  (86j:17016)}

\bibitem[Dzh85]{Dzhu:central-Zassenhaus}
\bysame, \emph{Central extensions of the {Z}assenhaus algebra and their
  irreducible representations}, Mat. Sb. (N.S.) \textbf{126(168)} (1985),
  no.~4, 473--489, 592. \MR{MR788083 (86j:17017)}

\bibitem[Far86a]{Far:associative}
Rolf Farnsteiner, \emph{The associative forms of the graded {C}artan type {L}ie
  algebras}, Trans. Amer. Math. Soc. \textbf{295} (1986), no.~1, 417--427.
  \MR{MR831207 (88f:17014)}

\bibitem[Far86b]{Far:central}
\bysame, \emph{Central extensions and invariant forms of graded {L}ie
  algebras}, Algebras Groups Geom. \textbf{3} (1986), no.~4, 431--455.
  \MR{MR901809 (88j:17014)}

\bibitem[Far87]{Far:dual}
\bysame, \emph{Dual space derivations and {$H\sp 2(L,F)$} of modular {L}ie
  algebras}, Canad. J. Math. \textbf{39} (1987), no.~5, 1078--1106.
  \MR{MR918588 (89a:17015)}

\bibitem[Jur99]{Ju:quotients}
G.~Jurman, \emph{Quotients of maximal class of thin {L}ie algebras. {T}he case
  of characteristic two}, Comm. Algebra \textbf{27} (1999), no.~12, 5749--5789.
  \MR{MR1726276 (2001a:17042b)}

\bibitem[Jur04]{Ju:Bi-Zassenhaus}
\bysame, \emph{A family of simple {L}ie algebras in characteristic two}, J.
  Algebra \textbf{271} (2004), no.~2, 454--481. \MR{MR2025538 (2005a:17015)}

\bibitem[Jur05]{Ju:maximal}
\bysame, \emph{Graded {L}ie algebras of maximal class. {III}}, J. Algebra
  \textbf{284} (2005), no.~2, 435--461. \MR{MR2114564 (2005k:17041)}

\bibitem[JY]{JuYo:quotients}
G.~Jurman and D.~S. Young, \emph{Quotients of maximal class of thin {Lie}
  algebras in characteristic two: errata and addendum}, preprint.

\bibitem[Kac74]{Kac}
V.~G. Kac, \emph{A description of the filtered {L}ie algebras with which graded
  {L}ie algebras of {C}artan type are associated}, Izv. Akad. Nauk SSSR Ser.
  Mat. \textbf{38} (1974), 800--834. \MR{MR0369452 (51 \#5685)}

\bibitem[KK94]{KoK1}
A.~I. Kostrikin and M.~I. Kuznetsov, \emph{On the structure of modular {L}ie
  algebras associated with analytic pro-{$p$}-groups}, Dokl. Akad. Nauk
  \textbf{339} (1994), no.~5, 591--593. \MR{MR1316518 (96a:17023)}

\bibitem[KK96]{KoK3}
\bysame, \emph{Finite-dimensional {L}ie algebras with a nonsingular
  derivation}, Algebra and analysis (Kazan, 1994), de Gruyter, Berlin, 1996,
  pp.~81--90. \MR{MR1465446 (98i:17023)}

\bibitem[KLGP97]{KL-GP}
G.~Klaas, C.~R. Leedham-Green, and W.~Plesken, \emph{Linear pro-{$p$}-groups of
  finite width}, Lecture Notes in Mathematics, vol. 1674, Springer-Verlag,
  Berlin, 1997. \MR{MR1483894 (98m:20028)}

\bibitem[Kos96]{Kos:beginnings}
Alexei~I. Kostrikin, \emph{The beginnings of modular {L}ie algebra theory},
  Group theory, algebra, and number theory (Saarbr\"ucken, 1993), de Gruyter,
  Berlin, 1996, pp.~13--52. \MR{MR1440203 (98e:17028)}

\bibitem[K{\v{S}}69]{KoSa}
A.~I. Kostrikin and I.~R. {\v{S}}afarevi{\v{c}}, \emph{Graded {L}ie algebras of
  finite characteristic}, Izv. Akad. Nauk SSSR Ser. Mat. \textbf{33} (1969),
  251--322. \MR{MR0252460 (40 \#5680)}

\bibitem[Kuz89]{Kuz:truncated}
M.~I. Kuznetsov, \emph{Truncated induced modules over transitive {L}ie algebras
  of characteristic {$p$}}, Izv. Akad. Nauk SSSR Ser. Mat. \textbf{53} (1989),
  no.~3, 557--589, 671. \MR{MR1013712 (90k:17037)}

\bibitem[LG94]{L-G}
C.~R. Leedham-Green, \emph{The structure of finite {$p$}-groups}, J. London
  Math. Soc. (2) \textbf{50} (1994), no.~1, 49--67. \MR{MR1277754 (95j:20022a)}

\bibitem[LGM84]{L-GMcK:maximal-classification}
C.~R. Leedham-Green and Susan McKay, \emph{On the classification of
  {$p$}-groups of maximal class}, Quart. J. Math. Oxford Ser. (2) \textbf{35}
  (1984), no.~139, 293--304. \MR{MR755666 (85m:20027)}

\bibitem[LGM02]{L-GMcKay}
C.~R. Leedham-Green and S.~McKay, \emph{The structure of groups of prime power
  order}, London Mathematical Society Monographs. New Series, vol.~27, Oxford
  University Press, Oxford, 2002, Oxford Science Publications. \MR{MR1918951
  (2003f:20028)}

\bibitem[LGN80]{L-GN}
C.~R. Leedham-Green and M.~F. Newman, \emph{Space groups and groups of
  prime-power order. {I}}, Arch. Math. (Basel) \textbf{35} (1980), no.~3,
  193--202. \MR{MR583590 (81m:20029)}

\bibitem[Luc78]{Lucas}
\`E. Lucas, \emph{Sur les congruences des nombres eul\'eriens et des
  coefficients diff\'erentiels des fonctions trigonom\'etriques, suivant un
  module premier}, Bull. Soc. Math. France \textbf{6} (1878), 49--54.

\bibitem[Mat99]{Mat:thin-groups}
Sandro Mattarei, \emph{Some thin pro-{$p$}-groups}, J. Algebra \textbf{220}
  (1999), no.~1, 56--72. \MR{MR1713453 (2000h:20049)}

\bibitem[Sel67]{Sel}
G.~B. Seligman, \emph{Modular {L}ie algebras}, Ergebnisse der Mathematik und
  ihrer Grenzgebiete, Band 40, Springer-Verlag New York, Inc., New York, 1967.
  \MR{MR0245627 (39 \#6933)}

\bibitem[SF88]{SF}
Helmut Strade and Rolf Farnsteiner, \emph{Modular {L}ie algebras and their
  representations}, Monographs and Textbooks in Pure and Applied Mathematics,
  vol. 116, Marcel Dekker Inc., New York, 1988. \MR{MR929682 (89h:17021)}

\bibitem[Sha94a]{Sha:max}
Aner Shalev, \emph{Simple {L}ie algebras and {L}ie algebras of maximal class},
  Arch. Math. (Basel) \textbf{63} (1994), no.~4, 297--301. \MR{MR1290602
  (95j:17025)}

\bibitem[Sha94b]{Sha:coclass}
\bysame, \emph{The structure of finite {$p$}-groups: effective proof of the
  coclass conjectures}, Invent. Math. \textbf{115} (1994), no.~2, 315--345.
  \MR{MR1258908 (95j:20022b)}

\bibitem[Skr90]{Skr:Hamiltonian}
S.~M. Skryabin, \emph{Classification of {H}amiltonian forms over algebras of
  divided powers}, Mat. Sb. \textbf{181} (1990), no.~1, 114--133. \MR{MR1048834
  (91j:17024)}

\bibitem[Skr98]{Skr:low}
Serge Skryabin, \emph{Toral rank one simple {L}ie algebras of low
  characteristics}, J. Algebra \textbf{200} (1998), no.~2, 650--700.
  \MR{MR1610680 (99b:17020)}

\bibitem[Str91]{Str:Block}
Helmut Strade, \emph{Representations of the {$(p\sp 2-1)$}-dimensional {L}ie
  algebras of {R}. {E}. {B}lock}, Canad. J. Math. \textbf{43} (1991), no.~3,
  580--616. \MR{MR1118011 (92g:17025)}

\bibitem[Str98]{Strade6}
H.~Strade, \emph{The classification of the simple modular {L}ie algebras. {VI}.
  {S}olving the final case}, Trans. Amer. Math. Soc. \textbf{350} (1998),
  no.~7, 2553--2628. \MR{MR1390047 (98j:17020)}

\bibitem[Str04]{Strade:book}
Helmut Strade, \emph{Simple {L}ie algebras over fields of positive
  characteristic. {I}}, de Gruyter Expositions in Mathematics, vol.~38, Walter
  de Gruyter \& Co., Berlin, 2004, Structure theory. \MR{MR2059133
  (2005c:17025)}

\bibitem[SW91]{SerWil}
Shirlei Serconek and Robert~Lee Wilson, \emph{Classification of forms of
  restricted simple {L}ie algebras of {C}artan type}, Comm. Algebra \textbf{19}
  (1991), no.~6, 1603--1628. \MR{MR1113954 (92e:17025)}

\bibitem[SZ92]{ShZe:finite-coclass}
A.~Shalev and E.~I. Zel{\cprime}manov, \emph{Pro-{$p$} groups of finite
  coclass}, Math. Proc. Cambridge Philos. Soc. \textbf{111} (1992), no.~3,
  417--421. \MR{MR1151320 (93e:20030)}

\bibitem[vdK73]{vdK}
W.~L.~J. van~der Kallen, \emph{Infinitesimally central extensions of
  {C}hevalley groups}, Springer-Verlag, Berlin, 1973, Lecture Notes in
  Mathematics, Vol. 356. \MR{MR0364484 (51 \#738)}

\bibitem[Wat91]{Wat}
William~C. Waterhouse, \emph{Automorphisms and twisted forms of generalized
  {W}itt {L}ie algebras}, Trans. Amer. Math. Soc. \textbf{327} (1991), no.~1,
  185--200. \MR{MR1038018 (91m:17025)}

\bibitem[Wil80]{Wil:type-S}
Robert~Lee Wilson, \emph{Simple {L}ie algebras of type {$S$}}, J. Algebra
  \textbf{62} (1980), no.~2, 292--298. \MR{MR563228 (81c:17019)}

\bibitem[Zus92]{Zus:current-central}
Paul Zusmanovich, \emph{Central extensions of current algebras}, Trans. Amer.
  Math. Soc. \textbf{334} (1992), no.~1, 143--152. \MR{MR1069751 (93a:17023)}

\end{thebibliography}

\end{document}